\documentclass[a4paper,11pt,reqno]{amsart} 
\usepackage[body={13cm,21cm}]{geometry}
\usepackage[initials]{amsrefs}

\usepackage{amssymb,amscd}
\usepackage[all]{xy}
\usepackage[mathscr]{eucal}
\usepackage{enumerate}

\numberwithin{equation}{section}
\theoremstyle{plain}
\newtheorem{thm}{Theorem}[section]
\newtheorem{cor}[thm]{Corollary}
\newtheorem{lem}[thm]{Lemma}
\newtheorem{prop}[thm]{Proposition}
\newtheorem*{problem*}{Problem}
\theoremstyle{definition}
\newtheorem{Def}[thm]{Definition}
\theoremstyle{remark}
\newtheorem{exa}[thm]{Example}
\newtheorem{exs}[thm]{Examples}
\newtheorem{rem}[thm]{Remark}
\newtheorem*{claim*}{Claim}
\newtheorem*{exa*}{Example}
\newtheorem*{rem*}{Remark}
\newtheorem*{rems*}{Remarks}

\providecommand{\C}[1]{\mathcal{#1}}
\providecommand{\D}[1]{\mathbb{#1}}

\providecommand{\F}[1]{\mathfrak{#1}}

\providecommand{\accol}[1]{\lbrace#1\rbrace}
\providecommand{\croch}[1]{\lbrack#1\rbrack}
\providecommand{\scal}[1]{\langle#1\rangle}
\providecommand{\quota}{/\!\!/}
\DeclareMathOperator{\Aut}{Aut}
\DeclareMathOperator{\Br}{Br}
\DeclareMathOperator{\Bru}{Br_{nr}}
\DeclareMathOperator{\car}{char}
\DeclareMathOperator{\centre}{centre}
\DeclareMathOperator{\codim}{codim}
\DeclareMathOperator{\Char}{\textit{\textsf{X}}}
\DeclareMathOperator{\disc}{disc}
\DeclareMathOperator{\Div}{Div}
\DeclareMathOperator{\divisor}{div}
\DeclareMathOperator{\End}{End}
\DeclareMathOperator{\Ext}{Ext}
\DeclareMathOperator{\Gal}{Gal}
\DeclareMathOperator{\GL}{GL}
\DeclareMathOperator{\Hom}{Hom}
\DeclareMathOperator{\id}{id}
\DeclareMathOperator{\im}{Im}
\DeclareMathOperator{\Ort}{O}
\DeclareMathOperator{\PGL}{PGL}
\DeclareMathOperator{\Pic}{Pic}
\DeclareMathOperator{\PSL}{PSL} 
\DeclareMathOperator{\rank}{rank}
\DeclareMathOperator{\Res}{Res}
\DeclareMathOperator{\SL}{SL}
\DeclareMathOperator{\SO}{SO}
\DeclareMathOperator{\Sp}{Sp}
\DeclareMathOperator{\Spec}{Spec}
\DeclareMathOperator{\Spin}{Spin}
\DeclareMathOperator{\tr}{Tr}
\begin{document}
\title[The rationality problem for fields of invariants]
        {The rationality problem for fields\\
        of invariants
        under linear algebraic groups\\
        (with special regards to the Brauer group)}
\author[J.-L.~Colliot-Th\'el\`ene]{Jean-Louis Colliot-Th\'el\`ene$^\ast$}
\author[J.-J.~Sansuc]{Jean-Jacques~Sansuc$^\dagger$}
\address{$^{\ast}$ Universit\'e Paris-Sud,
               91405 Orsay, France}
\address{$^{\dagger}$ IMJ, case 7012, Universit\'e Paris 7,
               2 place Jussieu, 75251 Paris, France}
\date{\footnotesize June 11, 2005}
\maketitle
\section*{Introduction}

Let $V$ be a vector space over an algebraically closed field $k$
of characteristic zero.
Let $G$ be a reductive subgroup of  $\GL(V)$. 
Assume that $G$ acts almost freely on $V$, i.e.\ that
for a general  $v$ in $V$, the stabilizer of $G$ at $v$
is trivial. Let 
$k(V)$ 
denote the field of rational functions on $V$.

\begin{problem*}
Is the field of invariants
$k(V)^G$ purely transcendental over $k$?
\end{problem*}

This is an old problem.  In these notes, it will be referred to as 
the \emph{purity problem}, or the \emph{rationality problem}, indifferently.
The question is open when $G$ is a connected group.
Even the case where  $G$ is the projective linear group 
$\PGL_n$ is not known, except for small values of $n$.
However, when
$G$ is finite, the question has a negative answer, as was shown by D.~Saltman in
1984. Saltman's paper was soon followed by a series of papers of F.~Bogomolov, then
by further papers of D.~Saltman.
Between 1986 and 1988 we ran seminars on their work, 
and in July 1988 one of us lectured on this topic at
the IX Escuela Latinoamericana de Matem\'aticas, held
in Santiago de Chile, and a set of notes was distributed. 
Over the years,  versions of these notes were circulated and used
as a complement to the original work of Saltman and Bogomolov.
We are
grateful to the organizers of the 2004 International Conference in Mumbai
for giving us the opportunity to 
publish a revised version of our text.
We have not tried to update  the notes systematically, but
we have added references to work done since 1988.
The rationality problem may also be raised over
a field $k$ which is not algebraically closed. 
The interested reader is referred to
\cites{Sw1,Vos,Sa1,Sw2,CS1,CS2,Mer}. 

Here is the list of sections. Sections \ref{sec.rationality} to \ref{sec.examples}
are devoted to general definitions and results,
whereas Sections \ref{sec.def.bru} to \ref{sec.homogeneous.spaces} concentrate on
the computation of the
unramified Brauer group.
\subsection*{Contents}
\begin{enumerate}[1.]
\item	
Rationality.
\item	
Quotients.
\item	
General techniques.
\item	
Examples.
\item	
The unramified Brauer group.
\item	
A general formula.
\item	 
Linear action of a finite group.
\item	
Multiplicative action of a finite group.
\item	
Homogeneous spaces.
\end{enumerate}

Section \ref{sec.rationality} introduce notions close to ``rationality'':
unirationality, stable rationality, retract rationality. Definitions and
properties are given in terms of function fields --- the algebraic side ---
and in terms of algebraic varieties --- the geometric side.

Section \ref{sec.quotients} deals with fields of invariants $K^G$ under a linear algebraic
group $G$ acting on a function field $K$. It  also deals with
``geometric models'' of such an action, i.e. a $G$-action  on
an integral $k$-variety $X$ such that $K$ coincides with the function field $k(X)$ of $X$
and $K^G $ coincides with the function field $k(X/G)$ 
of some ``good'' quotient $X/G$.  The most relevant concept
for these notes is that of an ``almost free'' action --- sometimes
also called ``generically free'' action.

Section \ref{sec.general.techniques} discusses two basic techniques: 
a ``slice" technique  and a lemma which has since then gone
into the literature as the ``no-name lemma".
The ``no-name lemma'' says that for ``almost-free'' linear $G$-actions
the (stable) purity question  depends only on $G$.
The ``slice method'' enables one to see the field of invariants
for an action of a given group $G$ as the field of invariants
for another action of a smaller group $H$.

Section \ref{sec.examples} considers linear actions of some particular groups for which
the (stable) rationality problem has a positive answer: 
some
finite groups ($\F{S}_n, \F{A}_5$), solvable groups, the so-called ``special groups" ($\GL_n,
\SL_n, \Sp_{2n}$), the 
orthogonal groups  $\Ort_n, \SO_n$. For the spinor groups $\Spin_n$ 
and for $\PGL_n$,
there
are results in low dimension.

The unramified Brauer group of a function field $K/k$  is defined
in Section \ref{sec.def.bru} and some of its
basic properties are given. A key property is that if $K$ is a purely transcendental  extension of $k$,
then the unramified Brauer group of $K$ (over $k$) is trivial.

In Section \ref{sec.bru.formula},  a  quite general formula of Bogomolov for
the unramified Brauer group of the field of invariants of
an (almost free) action of a group $G$ is given, first for 
$G$ finite then more generally for $G$ reductive.
The finite bicyclic subgroups of $G$ play a key r\^ole in this computation.

In Section \ref{sec.linear.action}, this formula is applied to the case of 
a linear action of a finite group.  It is further 
specialized to the case  of
 nilpotent groups of  class $2$, where it yields 
concrete examples of fields of linear invariants which are not
rational,
among which one finds Saltman's original example.

Section \ref{sec.multiplicative.action} discusses (twisted) multiplicative
invariants  of a finite group $G$, after Saltman and others. 
In favourable circumstances, one may view the associated field 
of $G$-invariants  as the field of invariants of a linear action of
a finite group $G'$ which is an extension of $G$.
This is used  to  produce  other types
of fields of linear invariants which are not
rational. 

Finally, Section \ref{sec.homogeneous.spaces} 
gives in some detail
Bogomolov's proof of the vanishing of the unramified Brauer group
of the function field 
of a quotient $G/H$, where $H$ is a connected subgroup
of a simply connected
group $G$. One first proves the vanishing of the unramified Brauer group
of the field $k(V)^H$ for an almost free linear representation $V$ of $H$.

Precise references to the papers of Saltman and Bogomolov are
given in the text. Some of our proofs differ from the original ones.

\subsection*{Known examples of non-rationality}

For the action  of a finite group $G$ on a purely
transcendental extension $K$ of $k=\D{C}$,
non-rational fields of invariants $K^G$
of the following types have been exhibited :
\begin{enumerate}[(i)]
\item
Linear action of a  nilpotent group of
order $p^7$ and class $2$.
\item
Multiplicative actions of  $G=(\D{Z}/p)^3$.
\item
Twisted multiplicative actions of $G=(\D{Z}/p)^2$.
\item
For an arbitrary action,  $G=\D{Z}/2$
(Clemens-Griffiths \cite{ClGr}, Artin-Mumford \cite{ArMu}).
\end{enumerate}
Over a non-algebraically closed field $k$, there exist examples of non-rational fields 
of invariants $K^G$
for an almost free linear action of the following kinds:
\begin{enumerate}[(i)]
\addtocounter{enumi}{4}
\item
$k=\D{Q}$, $G=\D{Z}/47$ (Swan \cite{Sw1}, Voskresenski{\u\i} \cite{Vos}).
\item
$k=\D{Q}$, $G=\D{Z}/8$ (Saltman \cite{Sa1}, Voskresenski{\u\i} \cite{Vos}).
\item
$k=\D{Q}$, $G$ a $k$-torus \cite{CS1}.
\item
$G$ a simply connected semisimple group (Merkurjev \cite{Mer})
\end{enumerate}
But for almost free linear
actions of connected linear algebraic groups over the complex field,
the rationality question for fields of invariants is open.

\section{Rationality}                        \label{sec.rationality}

Let $k$ be a
 field.  Let $\D{A}^n$, resp.\ $\D{P}^n$, denote the
$n$-dimensional affine
space, resp.\ projective
space, over $k$. The function field of an integral
$k$-variety $X$ is denoted $k(X)$. 
The set of rational $k$-points of $X$ is denoted $X(k)$.

One says that two integral $k$-varieties $X$ and $Y$ are
\emph{$k$-birationally equivalent} if the following equivalent conditions are
satisfied:
\begin{enumerate}[(i)]
\item
The function fields $k(X)$ and $k(Y)$ are isomorphic (over $k$).
\item
There exist non-empty Zariski open sets $U\subset X$ and $V\subset Y$
which are isomorphic over $k$.
\end{enumerate}
One says that two integral $k$-varieties $X$ and $Y$ are
\emph{stably $k$-birationally equivalent} if the following equivalent conditions are
satisfied:
\begin{enumerate}[(i)]
\item
For suitable integers $r,s$ and
independent variables $\accol{x_i}_{i=1}^r$ and $\accol{y_j}_{j=1}^s$,  the fields
$k(X)(x_1,\dots, x_r)$ and $k(Y)(y_1,\dots ,y_s)$ are $k$-isomorphic
(one then says that the fields $k(X)$ and $k(Y)$ are \emph{stably equivalent}).
\item
For suitable integers $r,s$,  the $k$-varieties $X\times_k\D{A}^r$
and $Y\times_k\D{A}^s$
are $k$-birationally equivalent.
\end{enumerate}
A $k$-variety $X$ is said to be 
\begin{itemize}
\item
\emph{$k$-rational} if it is integral and it is
$k$-birational to an affine space (one then says that $k(X)$ is pure over $k$),
\item
\emph{stably $k$-rational} if there exists an affine
space $\D{A}^n$ over $k$ such that
$X\times_{k} \D{A}^n$ is $k$-rational (one then says that $k(X)$ is 
stably pure over $k$),
\item
a \emph{direct factor of a $k$-rational variety} if
there exists an integral $k$-variety $Y$
such that $X\times_{k} Y$ is a $k$-rational variety.
\item
\emph{$k$-unirational} if it is integral and it
satisfies one of the
equivalent properties in the following lemma.
\end{itemize}

\begin{lem}                                                                \label{lem.unirational}
Let $X$ be an integral $k$-variety. 
The following conditions are
equivalent:
\begin{enumerate}[\rm(i)]
\item
The function field $k(X)$ of $X$ is a $k$-subfield of a pure
extension
$K$ of $k$.
\item
There exists a dominating $k$-morphism from a $k$-rational variety $Y$
to $X$.
\item
(Under the additional assumption that $k$ is infinite)
there exists a dominating $k$-morphism from a $k$-rational variety $Y$
to $X$, with $\dim Y = \dim X$.
\end{enumerate}
\end{lem}

\begin{proof}
We only need to show that (ii) implies (iii). 
Let $f\colon W \to X$ be a
dominating $k$-morphism, where $W \subset \D{A}^n$ is an open set
of affine space over $k$. By linear subspace of $W$ we shall mean
the non-empty trace on $W$ of a linear subspace of $\D{A}^n$.
Let us consider the closed, geometrically integral subvarieties
$Y\subset X$ over $k$ with the
following property: there exists a linear subspace $V\subset W$ with
$\dim V= \dim Y$ such that
the $k$-morphism $f$ restricts to a dominant $k$-morphism from $V$ to
$Y$.
Any  $k$-point in $f(W(k))$  is of this type. 
Let $Y$ be such a variety.  Assume $Y \neq X$.
Since $k$ is infinite and $f$ is dominant, there exists a $k$-point $P \in W(k)$
such that $M=f(P) \in  X(k)$  does not lie on $Y$.
Since $f$ is
defined at the generic
point of $V$, it is also defined at the generic point of the linear span
$L\subset W$ of $V$ and
$P$. Moreover, the closure $Y_1 \subset X$ of the image of $L$ under $f$ contains $Y$
and $M$, hence is of
dimension strictly bigger than $Y$. Iterating this procedure, we find that
$X$ is covered by a
$k$-linear space of dimension $\dim X$.
\end{proof}

Finally, following Saltman, one says that a $k$-variety $X$ is 
\begin{itemize}
\item
\emph{retract rational} (over $k$)  if it satisfies one
of the equivalent conditions in the following proposition (one then says that $k(X)$ is retract rational over $k$).
\end{itemize}

\begin{prop}[Saltman \cite{Sa2}*{Theorem 3.8}]          \label{prop.retract.rational}
Let $k$ be a field and  $X$ be an integral
$k$-variety. The
following conditions are equivalent:
\begin{enumerate}[\rm(i)]
\item
There exists a non-empty open set $U$ of $X$ such
that the identity morphism of $U$ factorizes through a Zariski open set $Y$
of an affine
space over $k$, \emph{i.e.}\ there are maps $U \to Y \to U$ whose composite is
identity on $U$.
\item
There exists a non-empty open set $V$ of $X$ such that for any
local
$k$-algebra $A$ with residue field $\kappa$, the natural map $V(A)\to
V(\kappa)$
is onto.
\end{enumerate}
\end{prop}

\begin{proof} 
The surjectivity of the natural map $V(A)\to V(\kappa)$ means that any
map
\[
\Spec\kappa\to V
\]
extends to $\Spec A$:
\[
\Spec\kappa\to\Spec A \to V.
\]
If $U$ is a Zariski open set of $V$, it is stable by generisation, hence
the surjectivity of
$V(A)\to V(\kappa)$ implies that  of $U(A)\to U(\kappa)$.

That (i) implies (ii) is now clear. For an affine space $Y=\D{A}^n$ the map
$Y(A)\to Y(\kappa)$ is the
natural reduction $A^n\to\kappa^n$ and  it is  surjective. Thus this is
also true for the open set
$Y$ of affine space, hence also for $U$.

In order to prove the converse assertion, we may assume that $X=V$ and that
$X$ is affine, say 
$X=\Spec B$ with $B=R/\F{p}$, with $\F{p}$ a prime ideal of
 the $k$-algebra $R=k[x_1,\dots,x_n]$ of
polynomials in $n$ variables
$x_1,\dots,x_n$. Let $K = k(X)$ be the function field of $X$. This is also
the residue class field of
the local ring $A= R_{\F{p}}$. The generic point $\eta$ of $X$ defines a
point in $X(K)$. The
assumption implies that this point comes from an element of $X(R_{\F{p}})$, i.e.\ we
have maps
\[
\eta\to\Spec R_{\F{p}} \to X
\]
whose composite is the natural inclusion of $\eta$ into $X$. The map
$\Spec R_{\F{p}} \to X$
factorizes through some $\Spec R_f$, where $R_f=R\croch{\frac{1}{f}}$, i.e.\ we have
an open subset
$V\subset\D{A}^n$ and maps
\[
\eta \to V\xrightarrow{\pi}X
\]
whose composite is the natural inclusion of $\eta$ into $X$. Then this
situation extends to an
open set $U$ of $X$, i.e.\ there exist a dense open set $U\subset X$ and
maps 
\[
U \to V\xrightarrow{\pi}X
\]
whose composite is the natural inclusion of $U$ into $X$. Setting
$Y=\pi^{-1}(U)$ gives rise to the announced factorization
\[
U \to Y \to U.\qedhere
\]
\end{proof}

\begin{rem}                          \label{rem.retract.rational}
Saltman's motivation for introducing the concept of retract
rationality was to try to
understand the relation between rationality and approximation properties
(in a more arithmetical
context). The criterion above, though nice, seems a priori to be of little
value: for instance it
seems impossible to prove that a smooth affine conic satisfies the lifting
property (ii) without a
priori proving that it is a rational curve! Nevertheless, in the particular
case of an almost free (see \S 2)
linear action of the projective linear group
$\PGL_p$ with  $p$ prime,
Saltman   used (ii) to
prove that the quotient is retract rational \cite{Sa2}*{Corollary 5.3} --- 
but a direct proof that this variety
satisfies (i) could later be given \cite{CS2}*{Corollary 9.13}.
\end{rem}

\begin{prop}                                             \label{prop.rationality.notions}
Let $X$ be an integral $k$-variety. 
In the list of properties:
\begin{enumerate}[\rm(i)]
\item
$X$ is $k$-rational,
\item
$X$ is stably $k$-rational,
\item
$X$ is a direct factor of a $k$-rational variety,
\item
$X$ is retract rational over $k$,
\item
$X$ is $k$-unirational,
\end{enumerate}
each property implies the following one.
\end{prop}

\begin{proof}  
Only the penultimate implication (iii)$\implies$(iv) 
requires a proof.
Let $Y$ be
an integral $k$-variety such that $X\times_{k} Y$ is rational. Let $U\subset
X\times_{k} Y$ be a non-empty open set
which is isomorphic to an open set of affine space. Let $(x_0,y_0)\in
U(k)$. Let $X_1\subset X$ be the non-empty open
set such that
$X_1\times\accol{y_0}=U\cap(X\times\accol{y_0})$. The open set $U_1= U\cap
(X_1\times Y)$ is
still isomorphic to an open set of affine space. 
Now the composite map $X_1\to U_1\to X_1$, where the first
map is given by $x\mapsto (x,y_0)$ and the second map is induced by
projection onto $X$, satisfies
the requirements of Proposition \ref{prop.retract.rational} (i).
\end{proof}

\begin{rems*}  
Over an algebraically closed field of characteristic zero, examples are known of:
\begin{enumerate}[(a)]
\item
unirational varieties which are not rational (Artin-Mumford \cite{ArMu},
Iskovs\-kikh-Manin \cite{IskMa}, Clemens-Griffiths \cite{ClGr});
\item
stably rational varieties which are not rational \cite{BCSS}.
\end{enumerate}
The Artin-Mumford method for proving non-rationality uses the Brauer group.
This method will be discussed and applied in \S\S 5\textendash 9 of these notes.
The method applies to varieties of arbitrary dimension. 
The associated invariant is insensitive to replacement of the variety $X$ by $X\times_{k} \D{A}^n$.
Such is not the case for the two other methods mentioned.

The Clemens-Griffiths method is quite specific to threefolds. 
That method, in a variant due to Mumford, was used to produce the
examples in (b).

The Iskovskikh-Manin method was originally developped
to prove the non-rationality of the general quartic threefold.
It has witnessed a strong development: birational rigidity and birational super-rigidity, 
work of Iskovs\-kikh \cite{Isk}, 
Pukhlikov \cites{Pu1,Pu2}, de Fernex-Ein-Musta\c{t}\u{a} \cite{dFEM}.

One should here also mention two further techniques for
disproving rationality. One is a
natural generalization of the Artin-Mumford technique:
it uses higher dimensional unramified cohomology (with
torsion coefficients). For this, see  \cites{CO1,O,Pe1,Pe2,C,Mer,Sa12,Sa13,Sa14}.
The other one, due to Koll\'ar \cites{Kol1, Kol2},  
uses reduction to positive characteristic.
\end{rems*}

\section{Quotients}                         \label{sec.quotients}

In this section $k$ denotes a field of characteristic zero and $\bar k$ an algebraic closure of $k$.
We denote by $\F{g}$ the Galois group of $\bar k$ over $k$.
Given a $k$-variety $X$ we denote $\bar X=X\times_k\bar k$. If $X$ is integral, resp.\ geometrically integral,
we denote by $k(X)$, resp.\ $\bar k(X)$, the function field of $X$, resp.~$\bar X$.

Let $G$ be a linear algebraic group over $k$ and let $X$ be
a geometrically integral $k$-variety with a $G$-action. Our interest will be in the field
$k(X)^G$, which by definition is the
fixed field of $\F{g}$ acting on the field of invariants $\bar k(X)^{G(\bar k)}$. 

\subsection{Elementary properties}  
                                      \label{ssec.elementary.properties}

The ring of regular functions on a $k$-variety $X$ is denoted $k\croch{X}$.

\begin{lem}                                           \label{lem.finite.invariant.functions}
Let $A$ be a domain, $K$ its field of fractions,
$G$ a finite group acting on $A$. Then the field
$K^G$ is the field of fractions of $A^G$.
\end{lem}

\begin{proof}
Let $f \in K^G$. Write $f=a/b$ with $a,b \in A$.
Let $e \in G$ be the identity element. Since $A$ is commutative, 
$\beta:=N(b)=\prod_{\sigma \in G}{}^ {\sigma}b\in A^G$.
Moreover, $\beta=bc$ with $c:=\prod_{\sigma \in G, \sigma \neq e}{}^ {\sigma}b\in A$. 
We have $f=\alpha/\beta$ with $\alpha:=ac\in A$.
Since $f$ and $\beta$ are $G$-invariants, $\alpha=f\beta$ is also $G$-invariant. 
Finally, $f=\alpha/\beta$ with
$\alpha,\beta\in A^G$.
\end{proof}

\begin{lem}                                                                              \label{lem.algebraic.group.invariant.functions}
Let $k=\bar k$. Let $A$ be an integral $k$-algebra of finite type,
equipped with an (algebraic) action of a linear algebraic
$k$-group $G$. Let $K$ be the field of fractions of $A$.
Assume that $A$ is a UFD and that
the identity component of
$G$ has no nontrivial character. 
Then $K^G$ is the field of fractions of $A^G$.
\end{lem}

See \citelist{\cite{Nag1}*{Theorem 4.1}\cite{PV}*{Lemma 3.2 and Theorem 3.3}}.

\begin{proof}
(i)
Let us first assume that $G$ is connected. For each prime ideal of height one $\F{p}$ of $A$,
fix a generator $f_{\F{p}} \in A$.
Let $f \in K^G$. Write
\[
f = uf_1^{r_1}\dots f_s^{r_s}
\]
with $u \in A^*$, the $f_i=f_{\F{p}_i}$ among the generators chosen above
and the $r_i\in\D{Z}$, $r_i\neq 0$. 
For any $\sigma\in G(k)$,  the identity ${}^{\sigma}f=f$ implies
\[
f=uf_1^{r_1}\dots f_s^{r_s}={}^{\sigma}u\cdot{}^{\sigma}f_1^{r_1}\dots{}^{\sigma}f_s^{r_s}.
\]
Then, by unique factorization, for each $i=1,\dots,s$ there exists a $j=\tau_{\sigma}(i)$ such that 
${}^{\sigma}f_i=\varepsilon_i(\sigma)f_j$, with $\varepsilon_i(\sigma) \in A^*$.
Since $G$ is connected, the induced homomorphism $\tau\colon G\to\F{S}_s$ is trivial.
For each $i$, we thus have ${}^{\sigma}f_i=\varepsilon_i(\sigma)f_i$. 
Let $\varepsilon_i\colon G(k) \to A^*$
be the induced homomorphism, and $\pi\colon A^*\to A^*/k^*$. 
Since $G$ is connected and $A^*/k^*$ is 
an abelian group 
of finite type (see \cite{CS1}), $\pi\circ\varepsilon_i$ is trivial. Hence $\varepsilon_i$ is a homomorphism
$\varepsilon_i\colon G(k) \to k^*$.
All data being algebraic, it is induced by a character $G \to\D{G}_m$, hence trivial by assumption.
Thus each $f_i$ is $G$-invariant.
This then implies $u\in A^G$, and $f$ lies in the fraction field of $A^G$.

(ii)
Let $G$ be arbitrary. Let $G^{\circ}$ be the identity component of $G$. 
Let $f\in K^G$. By (i), $f=a/b$ with $a,b\in A^{G^{\circ}}$.
The same arguments as in the proof of Lemma \ref{lem.finite.invariant.functions} 
for the finite group $G/G^{\circ}$ show that $f=\alpha/\beta$ with $\alpha\in A$ and $\beta\in A^G$, hence also
$\alpha\in A^G$.
\end{proof}

\begin{rems*}
(1) 
The above arguments also shows that $A^G$ is a UFD.

(2) The assumption that the identity component of
$G$ has no nontrivial character is a necessary one, 
as the example of the diagonal action of $\D{G}_m$ on $\D{A}_k^2$ shows.
\end{rems*}

\begin{prop}                                               \label{prop.invariant.functions}
Let $k=\bar k$. Let $G$ be a linear algebraic group and $X$ a 
factorial affine variety with a $G$-action.
\begin{enumerate}[\rm(i)]
\item
If $k[X]^*=k^*$, there exists a $G$-invariant affine open set
$U\subset X$ such that $k(X)^G$ is the field of fractions of $k\croch{U}^G$.
\item
If  $G$ is
finite, or if $G$ has no nontrivial character $G\xrightarrow{\chi}\D{G}_m$, then $k(X)^G$
is the field of fractions of $k[X]^G$.
\end{enumerate}
\end{prop}

See also \cite{PV}*{Theorem 3.3}.

\begin{proof}
(i)
Let $k(X)^G=f(f_1,\dots,f_s)$. Let $Z$ be the union of the supports of the divisors of the functions $f_i$.
Since $X$ is factorial affine, there exists $g\in k\croch{X}$ whose divisor is $Z$. 
The complement $U$ of $Z$ is therefore an affine open set with $k\croch{U}=k\croch{X}\croch{1/g}$.
The divisor $Z$ being $G$-invariant, $U$ is $G$-invariant. Since $f_1,\dots,f_s$ belong to $k\croch{U}^G$,
this implies $k(X)^G)$ is the field of fractions of $k\croch{U}^G$.

Assertion (ii) is just a rephrasing of  Lemmas \ref{lem.finite.invariant.functions} 
and \ref{lem.algebraic.group.invariant.functions}.
\end{proof}

In order to describe a field of invariants $k(X)^G$, it will often be convenient to realize it
as the
function field of a ``quotient'' variety $Y$ of $X$ by $G$. 

In the literature, several integral
varieties $Y$ equipped with a morphism $p\colon X\to Y$ such that $p(gx)=p(x)$
for $g\in G$ and $x\in X$
go under the name of quotient varieties. However, our main interest is in
such quotients with the
additional property $k(Y)=  k(X)^G$. 

\begin{prop}                                   \label{prop.generic.quotient}
Let $\pi\colon  X\to Y$ be a dominant
morphism of geometrically
integral $k$-varieties. Let the algebraic $k$-group $G$ act on $X$.  Assume that
for $x, y\in X(\bar k)$, the following two conditions are equivalent:
\begin{enumerate}[\rm(i)]
\item
$\pi(x)=\pi(y)$;
\item
there exists $g\in G(\bar k)$ such that $gx=y$.
\end{enumerate}
Then $k(Y)$ may be identified with the field
$k(X)^G$ of invariants:
\[
k(Y)=k(X)^G.
\]
\end{prop}

\begin{proof}
For the proof, we may assume $k=\bar k$. In this case the proposition is an immediate consequence of the
following elementary
proposition --- whose proof uses the characteristic zero hypothesis.
\end{proof}

\begin{prop}[see \cite{Kr}*{AI.3.7 Satz 2}]    \label{prop.dominant.quotient}
Let $k=\bar k$. Let $\phi\colon X\to Y$ be a
dominant morphism of integral $k$-varieties. If $f\in k(X)$ is a rational
function which is
constant on the fibres of $\phi$, then $f$ may be identified with a
rational function on $Y$.
\end{prop}

\subsection{Algebraic quotients, geometric 
           quotients,  torsors}    \label{ssec.agt.quotients}

\begin{Def}                                                               \label{def.algebraic.quotient}
If $G$ is a reductive $k$-group and $X=\Spec A$
is an affine
$k$-variety with a $G$-action, the \emph{algebraic quotient} $X\quota G$ of $X$ by $G$
is the
affine scheme $\Spec A^G$. This actually is a $k$-variety, since $G$ reductive
implies (see \citelist{\cite{Kr}*{II.3.2}\cite{MuFo}*{Chap.~1, \S2}})
that
$A^G$ is a finitely generated $k$-algebra.
\end{Def}

\begin{Def}                                                               \label{def.geometric.quotient}
Let $X$ be an algebraic $k$-variety equipped 
with an action of an
algebraic $k$-group $G$.  A \emph{geometric quotient} of $X$ by $G$
is a $k$-variety $Y$ equipped with a $k$-morphism $\phi\colon X\to Y$ such that:
\begin{enumerate}[(i)]
\item
$\phi$ is open, constant on $G$-orbits and induces a bijection of $X(\bar k)/G(\bar k)$ with $Y(\bar k)$.
\item
For any open subset $V\subset Y$, the natural morphism $k\croch{V}\to k\croch{\phi^{-1}(V)}^G$
is an isomorphism.
\end{enumerate}
If such a variety $Y$ exists, it is unique up to unique isomorphism. It will then be denoted  $Y=X/G$.
\end{Def}

By a classical theorem of Rosenlicht \cite{R2}, given $X$ an integral algebraic variety equipped with an action of a linear algebraic group $G$, there exists a nonempty $G$-stable open set $U\subset X$ admitting a geometric quotient
$U\to U/G$.

\begin{rem*}                    
A geometric quotient $(Y,\phi)$ of $X$ by $G$ has the following properties:
\begin{enumerate}[(a)]
\item
The map $\phi\colon X(\bar k)\to Y(\bar k)$ is onto, all $G(\bar k)$-orbits on $X(\bar k)$ are closed and they coincide
with the fibres of $\phi$ over $\bar k$-points. 
\item
Proposition \ref{prop.generic.quotient} implies $k(Y)=k(X)^G$.
\item
The map $\phi\colon X\to Y$ has the obvious universal mapping property \cite{Bor}*{p.~172}. 
This accounts for the unicity statement above.
\item
When $X$ is affine and $G$ is reductive, the algebraic quotient $X\quota G$ need not be a geometric quotient 
(consider the diagonal action of $\D{G}_m$ on $\D{A}^2$).
\end{enumerate}
\end{rem*}

\begin{prop}                                                               \label{prop.geometric.quotient}
Let $k=\bar k$. Let $X$ be a normal algebraic $k$-variety equipped with an action of an
algebraic $k$-group $G$. 
Then a $k$-morphism $\phi\colon X\to Y$ is a geometric quotient if and only if
\begin{enumerate}[\rm(1)]
\item
$\phi$ is constant on $G$-orbits.
\item
$\phi$ induces a bijection $X(k)/G(k)\xrightarrow{\cong}Y(k)$.
\item
$Y$ is normal.
\end{enumerate}
\end{prop}

\begin{proof}                    
The proof uses the $\car k=0$ hypothesis and  Zariski's main theorem, see \cite{MuFo}*{Prop. 0.2, p.~7}.
\end{proof}

\begin{Def}                                                                                          \label{def.torsor}
Let $G$ be a linear algebraic $k$-group. 
A $k$-variety $X$ equipped with a $G$-action and a faithfully flat
$k$-morphism $X\to Y$ is a $G$-\emph{torsor} --- \emph{principal homogeneous
space} under $G$ --- over $Y$ if
the map $(g,x)\mapsto (gx,x)$ induces an isomorphism 
$G\times_k X\cong X\times_YX$. 
\end{Def}

\begin{rem}                    \label{rem.affine.torsor}
If $X$ is affine,
then $Y$ is affine and
$Y=X\quota G$.
\end{rem}                    

\begin{rem}                    \label{rem.geometric.quotient.torsor}
If $X$ is a $G$-torsor over $Y$ the map $X\to Y$
makes $Y$ into a geometric quotient $X/G$. 
Hence $k(Y)=k(X)^G$.
\end{rem}                    

\begin{rem}                    \label{rem.homogeneous.space.torsor}
Let $G$ be a $k$-algebraic group and $H$ a subgroup. 
The geometric quotient $G/H$ exists, and $G\to G/H$ makes $G$ into an $H$-torsor over $G/H$.
Hence $k(G/H)=k(G)^H$.
\end{rem}                    

There is a vast literature on the various possible notions of quotient. We refer the
reader to \cites{Bia1, Bia2, DiR, Kr,PV}.

\subsection{Almost free actions}  \label{ssec.almost.free}

\begin{Def}                                                                       \label{def.almost.free}
An action of an algebraic $k$-group $G$ on a geometrically integral
$k$-variety $X$
will be called an \emph{almost free action} 
if for all $x\in U(\bar k)$ in a non-empty open set $U$ of $X$ the stabilizer $G_x\subset G(\bar k)$ is trivial.
\end{Def}

\begin{thm}[Luna \cite{Lu}]                                         \label{thm.luna.torsor}
Let $X$ be an affine $k$-variety
and let the reductive $k$-group $G$ act on $X$. If each stabilizer $G_x$, $x \in X(\bar k)$,
is trivial,
then $X$ is a $G$-torsor over $X\quota G$.
Hence, $X/G$ exists, and $X/G=X\quota G$.
\end{thm}

\begin{rem*}                                                                
The hypothesis ``$X$ affine'' cannot be ignored
(see \cite{MuFo}*{p.~11}).
\end{rem*}

\begin{thm}[Luna \cite{Lu}]                                                         \label{thm.luna.slice}
Let $k=\bar k$. Let the reductive group $G$ act
upon the
affine variety $X$.  Let $\pi\colon X\to X\quota G$ be the natural projection. If
there
exists a point $x\in X(k)$ with trivial stabilizer and closed orbit, then
there exists an
affine open set $V\subset X\quota G$ which contains $\pi(x)$ and is such that the
projection $\pi^{-1}(U)\to U$ makes $\pi^{-1}(U)$ into a $G$-torsor.
\end{thm}

\begin{proof} 
This is an immediate consequence of Luna's ``slice
theorem'' \cite{Lu}.
\end{proof}

\begin{cor}                                                                                \label{cor.luna.slice}
Let $k=\bar k$. Let the reductive group $G$ act
upon the
 variety $X$. 
 The following
conditions are equivalent:
\begin{enumerate}[\rm(i)]
\item
There exist a $G$-invariant non-empty open set $U\subset X$ and a variety $Y$
with a
morphism $U\to Y$ such that $U$ is a $G$-torsor over $Y$.
\item
There exists a $G$-invariant non-empty affine open set $U\subset X$ such 
that the
projection map $U\to U\quota G$ makes $U$ into a $G$-torsor over $U\quota G$.
\item
There exist an \emph{affine} $G$-invariant open set $U$ and a point
$x\in U(k)$ such
that the stabilizer $G_x$ is trivial and the orbit of $x$ in $U$ is
closed.
\end{enumerate}
\end{cor}

\begin{proof} 
To show that (i)$\implies$(ii), choose an affine non-empty open set
$V\subset Y$
and take its inverse image $\pi^{-1}(V)$ under $\pi\colon U\to Y$. Conversely, (ii)$\implies$(i) is obvious.
It is also obvious that (ii)$\implies$(iii). Finally, (iii)$\implies$(ii) is a consequence of Theorem \ref{thm.luna.slice}.
\end{proof}

\begin{thm}[Popov \cite{Po1}, 
Luna-Vust \cite{LuVu}, see \cite{MuFo}*{appendix p.~154}]         \label{thm.nagata.semisimple}
Let $k=\bar k$. Let the semisimple group $G$ act on
the factorial affine variety
$X$.
The following conditions are equivalent:
\begin{enumerate}[\rm(i)]
\item
There exists a $G$-invariant non-empty open set $U\subset X$ such that for
each
$x\in U(k)$, the stabilizer $G_x$ is reductive.
\item
There exists a $G$-invariant non-empty open set $V\subset X$ such that
for each $x\in V(k)$, the orbit $G\cdot x$ of $x$ is closed in $X$.
\end{enumerate}
\end{thm}

As the diagonal action of the multiplicative group $\D{G}_m$ on affine space
reveals, 
Theorem \ref{thm.nagata.semisimple} does not extend to reductive groups. We thank M.~van den Bergh for
showing us
how to prove the next theorem for arbitrary reductive groups.

\begin{thm}                                                                    \label{thm.vdbergh.reductif}
Let $k=\bar k$. Let $X$ be a factorial affine variety,
and let
the reductive group $G$ act on $X$. Assume that the open 
set consisting of
points with
finite stabilizer is non-empty.

Then there exist non-empty $G$-invariant open sets $U_1\subset U\subset X$
such that $U$ is affine and $G$-orbits of points of $U_1$ are closed in
$U$.
\end{thm}

\begin{proof}
Let $G^{\circ}\subset G$ be the connected component of
identity in $G$. Assume we
have found suitable $U_1\subset U$ for $G^{\circ}$. Let $V$ be the
intersection of all
$\sigma U$, for $\sigma$ running through a (finite) system of
representatives in $G$ for
$G/G^{\circ}$. The open set $V$ is affine and $G$-stable. Let $V_1=U_1\cap
V$. The $G$-orbits of
points of $V_1$ are finite unions of orbits of $G^{\circ}$, hence are closed
in $U$, hence also in
$V$. We may therefore assume that $G$ is connected.

Let $T=Z(G)^{\circ}$ be the torus which is the connected centre of $G$, and
let $G_0=\croch{G,G}$ be the
derived subgroup of $G$. It is a connected semisimple group. Then $G$ is
the almost direct product
of $T$ and $G_0$, i.e.\ it is the quotient of $T\times G_0$ by the finite
group $ =T\cap G_0$. Let
$f\colon X\to X\quota G_0$ be the algebraic quotient of $X$ by $G_0$. The torus $T$
acts on $X\quota G_0$ and the map
$f$ is $T$-equivariant. Any point of $X$ with finite stabilizer in $G$
projects down to a point of
$X\quota G_0$ with finite stabilizer in $T$. Hence there exists a non-empty
$T$-stable open set
$W_0\subset X\quota G_0$ consisting of points with finite stabilizer. All orbits
of $T$ on $W_0$ are
closed (if an orbit of $T$ on $W_0$ were not closed, its closure would
contain a smaller dimensional
orbit, hence points with positive dimensional stabilizer).

Since $T$ is a torus, a theorem of Sumihiro \cites{Su} asserts that there exists
a non-empty affine
$T$-invariant open set $U_0\subset W_0$. Clearly, $T$-orbits are closed in
$U_0$. Let
$U=f^{-1}(U_0)$. Since $f$ is affine, this is a $G$-invariant affine open
set of $X$. Since $X$ is
factorial, so is $U$. By assumption, $X$ and hence $U$ contains a non-empty
$G$-stable open set
consisting of points with finite stabilizer. Thus according to Theorem \ref{thm.nagata.semisimple}
there exists a non-empty
$G_0$-stable open set $U_1$ of $U$ consisting of points with closed
$G_0$-orbit in $U$ and finite
stabilizer. For any point $x\in U_1(k)$, $G_0x$ is of maximal dimension in $f^{-1}(f(x))$.
Since it is
closed, it coincides with $f^{-1}(f(x))$. Now
\[
Gx=TG_0x=T(f^{-1}(f(x)))= f^{-1}(Tf(x)))
\]
is the inverse image of the closed set $Tf(x)\subset U_1$, hence is
closed.
\end{proof}

\begin{cor}                                                 \label{cor.vdbergh.reductif}
Let $k=\bar k$. Let $X$ be a factorial affine variety,
and let the reductive group $G$ act on $X$. Suppose $G$ acts almost freely.
Then  there exists a $G$-invariant non-empty 
affine open set $U\subset X$ such that the
projection map $U\to U\quota G$ makes $U$ into a $G$-torsor over $U\quota G$.
\end{cor}

\begin{proof} 
According to Theorem \ref{thm.vdbergh.reductif}, there exist non-empty
$G$-invariant open sets
$U_1\subset U\subset X$ such that $U$ is affine and $U_1$ consists of
points $x$ with trivial
stabilizer and closed orbit in $U$. It only remains to apply Corollary
\ref{cor.luna.slice}, (iii)$\implies$(i).
\end{proof}

\section{General techniques}
                                           \label{sec.general.techniques}

In this section, $k$ denotes a field of characteristic zero and $\bar k$ an algebraic closure of $k$.

Several methods have been used to prove the rationality of some quotient
spaces, in particular by Katsylo. A good survey of
these methods is given by Dolgachev \cite{Do}, to whom we refer for many
examples from the theory of
moduli spaces. In this section, these general methods will be reviewed; an
attempt will be made to
give complete proofs of some of the basic lemmas.

\subsection{The slice method}  \label{ssec.slice.method}

The following very useful fact, which is known under the name ``slice
method'' goes back at least to C.~S.~Seshadri 
\cite{Sesh}. See M.~Nagata \cite{Nag2}*{Lemma 1 (``Lemma of Seshadri''), p.~37} 
and also V.~L.~Popov \cite{Po4}. 
The following version of Seshadri's lemma comes from \cite{GaVi}.

\begin{thm}[Slice lemma]                                                         \label{thm.slice.lemma}
Let $G$ be an algebraic group over $k$ and  $X$
a geometrically integral
$k$-variety with a $G$-action. Let $Y\subset X$ be a closed geometrically integral
sub-$k$-variety, and let
$H\subset G$ be the normalizer of $Y$. Assume:
\begin{enumerate}[\rm(i)]
\item
The closure of $G\cdot Y$ is the whole of $X$.
\item
There exists a non-empty $H$-stable open subset $Y_0$ of $Y$
such that
if $g\in G(\bar k)$ and $y\in Y_0(\bar k)$ satisfy $g\cdot y\in Y_0(\bar k)$, then there exists
$h\in H(\bar k)$
such that $h\cdot y=g\cdot y$.
\end{enumerate}
Then there is a natural $k$-isomorphism of fields $k(X)^G\cong k(Y)^H$.
\end{thm}

\begin{proof} 
We immediately reduce to the case $k=\bar k$.
Let $f\in k(X)^G$, $f\neq 0$. Since $G\cdot Y$ is dense
in $X$, the subvariety $Y$ is not contained
in the support of $\divisor f$. Hence $f|_Y\in k(Y)$ and the restriction
defines an injection $k(X)^G
\subset k(Y)^H$.

Conversely, let $f\in k(Y)^H=k(Y_0)^H$. Let $\pi\colon G\times Y_0\to X$ be the
natural morphism. The
rational function $\phi$ defined on $G\times Y_0$ by $\phi (g,y) = f(y)$ is
constant on the fibres of
$\pi$ (use assumption (ii) and the $H$-invariance of $f$). Since $\pi$ is a
dominant morphism of integral
varieties, $\phi$ may be identified --- see Proposition \ref{prop.dominant.quotient} ---
with a rational function $\tilde f$ on
$X$, and $\tilde f|_Y =
f$. Hence $k(X)^G = k(Y)^H$.
\end{proof}

\begin{rem}                                              \label{rem.popov.chevalley.slice}
We refer to  \cite{Po5} for more historical 
remarks regarding this theorem, which may be consider as a birational analogue of the Chevalley restriction 
theorem \cite{LuR}.
\end{rem}

For other reduction techniques, see \cites{Kor2,Kor3}.

\subsection{The no-name lemma 
                      (\cite{BoKat,Do})}  \label{ssec.noname.lemma}

In quite a few circumstances, the method to be described below enables one
to deduce stable
rationality and in some cases rationality of some fields of invariants once
it is known for other
fields. This method has been independently discovered by several people, hence the 
denomination ``no-name lemma'' given by
Dolgachev \cite{Do} to one of the corollaries below.

We first give statements and proofs for finite groups actions. Then we give an
independent proof for the
general case of reductive group actions.

\begin{thm}[Speiser's lemma]    \label{thm.noname.finite}
Let $G$ be a
finite
group. Let $V$ and $W$ be two faithful finite-dimensional linear representations of $G$ over $k$. Then
the fields $k(V)^G$ and
$k(W)^G$ are stably equivalent.
\end{thm}

This means that there exist independent variables $x_1,\ldots,x_r$ and
$y_1,\ldots,y_s$ such that
\[
k(V)^G(x_1,\ldots,x_r)\cong k(W)^G(y_1,\ldots,y_s).
\]

\begin{cor}                                                         \label{cor.universal.stable.purity}
Let $V$ be a faithful finite-dimensional linear representation 
of a finite group $G$ over $k$. 
The stable purity of the field $k(V)^G$ depends only on~$G$.
\end{cor}

We have the well-known lemma:

\begin{lem}                                                               \label{lem.semilinear.action}
Let $K/k$ be a finite Galois extension. Let
$G$ be its Galois group.
Then for any finite dimensional $K$-linear space $E$ and any semi-linear action of $G$ on $E$
with
respect to the extension $K/k$
\[
E = E^G\otimes_kK.
\]
In particular, $K(E)^G/k = k(E^G)/k$ is pure.
\end{lem}

\begin{proof}
Let us prove the last assertion. If $E = E^G\otimes_kK$, we are
reduced to the case $E =
K^n$ with the natural action of $G$. Hence $K(E) = K(t_1,\ldots,t_n)$ with
trivial action of $G$
on the $t_i$'s. So $K(E)^G/k = k(t_1,\ldots,t_n)/k$. 

For the first assertion,
let
$\{\omega_i\}_{i=1,\ldots,d}$ be a linear basis of $K/k$ and
$G=\{\sigma_i\}_{i=1,\ldots,d}$ with
$\sigma_1=\id_K$. Let $v \in E$. Consider, for every $i=1,\ldots,d$,
\[
w_i:= \sum_j \sigma_j(\omega_iv) \in E^G.
\]
The group $G$ acting semi-linearly on $E$ with respect to $K/k$, we have
\[
w_i = \sum_j \sigma_j(\omega_i)\sigma_j(v).
\]
The matrix $(\sigma_j(\omega_i))$ being invertible, one can express $v =
\sigma_1(v)$ as a linear
combination of the $w_j$'s with coefficients in $K$, i.e.\ $v\in
E^G\otimes_kK$. Thus the natural map
\[
E^G\otimes_kK\to E
\]
is onto. It remains to show that  this map is injective. Let $\sum_i\lambda_iv_i=0$ be a non-trivial $K$-linear relation between elements $v_i\in E^G$. There exist $i_0$ and $\alpha\in K$ with $\tr_{K/k}(\alpha\lambda_{i_0})\neq 0$.
Then $\sum_i\tr_{K/k}(\alpha\lambda_i)v_i=0$ is a non-trivial $k$-linear relation between the $v_i$'s.
\end{proof}

\begin{proof}[Proof of Theorem \ref{thm.noname.finite}] 
By hypothesis, the linear
representation of $G$ on $V$ is
faithful. As a consequence the extension $k(V)/k(V)^G$ is Galois of group
$G$. Then we can
apply Lemma \ref{lem.semilinear.action} to the following data: the extension $k(V)/k(V)^G$ and 
the $k(V)$-linear space $E=W_{k(V)}$. 
Lemma \ref{lem.semilinear.action}
ensures the purity of the extension
\[
k(V\oplus W)^G/k(V)^G = (k(V)(W))^G/k(V)^G = K(W)^G/k.
\]
By exchange of $V$ and $W$ we also obtain the purity of $k(V\oplus
W)^G/k(W)^G$. Hence
$k(V\oplus W)^G$ is a common pure extension of $k(V)^G$ and of $k(W)^G$.
\end{proof}

We now give statements and proofs for the general case.

\begin{thm}                                                       \label{thm.noname.general}
Let $G$ be a reductive group 
over $k$.
Let $Y =
\Spec A$ be an
affine $k$-variety with a $G$-action. Assume that all  geometric stabilizers at closed points of 
$X$ are trivial.
Let $X\to Y$ be a vector bundle over $Y$ equipped with an
equivariant $G$-action. Then:
\begin{enumerate}[\rm(a)]
\item
$X\quota G$
is a vector bundle over $Y\quota G$.
\item
If $Y$ and hence $X$ are integral, the field
$k(X)^G$ is pure over $k(Y)^G$.
\end{enumerate}
\end{thm}

\begin{rem}                                                      \label{rem.noname.general}
In the paper \cite{BoKat} by 
Bogomolov and Katsylo, no
condition is imposed upon
the linear group $G$. Also, the base
variety $Y$ is not assumed to be affine.
Although this last extension is perfectly legitimate, one should be aware
that there are actions of
$\SL(2)$ on non-affine varieties such that:
\begin{enumerate}[(i)]
\item
All stabilizers are trivial.
\item
There is a geometric quotient.
\item
The action is not free (Mumford, \cite{MuFo}*{p.~11}).
\end{enumerate}
For a proof of part (b) for arbitrary $G$ and $Y$, see \cite{CGR}*{\S 4.3}.
\end{rem}

\begin{proof}[Proof of Theorem \ref{thm.noname.general}]
Let $M$ be the group of global
sections of the vector bundle
$X$ over $Y=\Spec A$. This is a finitely generated projective $A$-module, and
$X = V(M) = \Spec
\C{S}_A(M)$, where $\C{S}_A(M)$ denotes the symmetric algebra of $M$ over $A$.
Let $B=A^G$. 
Our assumptions imply that $Y\to Y\quota G$ is a $G$-torsor, see Theorem \ref{thm.luna.torsor}. 
Then, $A$ is a faithfully flat extension of $B$. By assumption,
there is a
$G$-equivariant action on the $A$-module $M$. This assumption may be
translated into descent data on
the projective $A$-module $M$. Since such descent data are effective, 
there
exists a projective
$B$-module $N$ such that $A\otimes_BN\cong M$, and this $B$-module $N$ is none other
but $M^G$. One then has
$\C{S}_A(M)^G=\C{S}_B(N)$. Since $X\quota G= \Spec\C{S}_A(M)^G$ and $Y\quota G= \Spec B$, it
follows that $X\quota G$ is a vector
bundle over $Y\quota G$, which is the first assertion. Now note that since all
geometric stabilizers of $G$ at a closed point of $Y$ are trivial, 
certainly the same holds for $X$, and $X$ is a
principal homogeneous space over
$X\quota G$, just as $Y$ is over $Y\quota G$. If $Y$ and hence $X$ are integral, so are
$X\quota G$ and $Y\quota G$, and
$k(X)^G=k(X\quota G)$ is purely transcendental over $k(Y\quota G)=k(Y)^G$.
\end{proof}

\begin{cor}                                                          \label{cor.noname.general}
Let $G$ be a reductive group over $k$, let $V$
be a finite-dimensional $k$-vector
space with a linear $G$-action and let $X$ be a factorial affine $k$-variety with an almost
free $G$-action. Then the field $k(X\times V)^G$ is pure over the field $k(X)^G$. In
particular,
if $k(X)^G$ is pure, so is $k(X\times V)^G$.
\end{cor}

\begin{proof} 
By Corollary \ref{cor.vdbergh.reductif} 
there exists a non-empty
$G$-invariant affine open set
$U\subset X$ such that $U\to U\quota G$ is a principal homogeneous space under
$G$, and $k(U\quota G)=k(X)^G$.
The natural projection from $Z_1=U\times V$ to $U$ makes $Z_1$ into a
(trivial) vector bundle over
the affine variety $U$. This vector bundle is equipped with the diagonal
$G$-action. According to  Theorem \ref{thm.noname.general}, 
$Z_1\quota G$ is a vector bundle over $U\quota G$. Since such a bundle is locally
trivial in the Zariski
topology, the function field $k(Z_1\quota G)$ is purely transcendental over
$k(U\quota G)=k(X)^G$. On the other
hand, for the diagonal action of $G$ on the affine variety $Z_1=U\times Y$
all stabilizers are
trivial, hence by Luna's Theorem \ref{thm.luna.torsor}, the map $Z_1\to Z_1\quota G$ makes $Z_1$ into a
principal homogeneous
space over $Z_1$ under $G$. In particular, $k(Z_1\quota G)=k(Z_1)^G=k(X\times
V)^G$. Thus $k(X\times V)^G$
is purely transcendental over $k(X)^G$.
\end{proof}

\begin{cor}[No-name lemma]                               \label{cor.noname.lemma}
Let $G$ be a
reductive group over $k$ and let
$V$ and $W$ be two finite-dimensional $k$-vector spaces with almost free linear $G$-actions. Then the
fields
$k(V)^G$ and $k(W)^G$ are stably equivalent. If one of them is pure, the
other one is
stably pure.
\end{cor}

\begin{proof} 
The diagonal action of $G$ on $V\times W$ gives a common extension:
\[
\xymatrix{
                                  &k(V\times W)^G&\\
k(V)^G\ar@{-}[ur]_{\text{pure}}&                             &k(W)^G,\ar@{-}[ul]^{\text{pure}}  
}
\]
where purity follows from Corollary \ref{cor.noname.general}.
\end{proof}

\begin{cor}                                                          \label{cor.appli.noname.lemma}
Let $G$ be a
reductive $k$-group and let $V$, $V_1$ and $V_2$ be finite-dimensional $k$-vector
spaces with linear 
$G$-actions. Assume that the action on $V$ is almost free. If the action of
$G$ on
$V_1$ is almost free and $\dim V_2 \geq \dim V$, then $k(V_1\oplus V_2)^G$
is pure
over $k(V)^G$. If moreover $k(V)^G$ is pure, so is $k(V_1\oplus V_2)^G$.
\end{cor}

\begin{proof}
Because the action of $G$ on $V_1$ is almost free,
the same argument as in the
previous proof shows that $k(V_1\oplus V)^G$ is pure over $k(V_1)^G$, of
transcendence degree $\dim
V$. Similarly, $k(V_1\oplus V_2)^G$ is pure over $k(V_1)^G$, of
transcendence degree $\dim V_2$. Thus
$k(V_1\oplus V_2)^G$ is pure over $k(V_1\oplus V)^G$. Now since the action
of $G$ on $V$ is almost
free, $k(V_1\oplus V)^G$ is pure over $k(V)^G$.
\end{proof}

\begin{cor}                                                  \label{cor.universal.almost.free.action}
Let the reductive $k$-group $G$ act almost
freely and linearly on the finite-dimensional $k$-vector
space $V$. Let $G\subset \GL_n$ be a closed embedding of $G$ as a subgroup of $\GL_n$ 
for some
integer
$n$. Then the field $k(V)^G$ is stably equivalent to the function field
$k(\GL_n/G)$
of the homogeneous space $\GL_n/G$.
\end{cor}

\begin{proof}
Let $\C{M}_n$ denote be the set of $n\times n$ matrices. Via the natural
(left) multiplication,
there is a linear $G$-action on the vector space $\C{M}_n$. On the $G$-stable
affine open set
$U=\GL_n\subset \C{M}_n$, the group $G$ acts with trivial stabilizers. 
Thus the action of $G$ on $\C{M}_n$ is almost free, and Corollary \ref{cor.noname.lemma}
implies
that $k(V)^G$ is stably
equivalent to the function field $k(\GL_n)^G$
which, by Remark \ref{rem.homogeneous.space.torsor}, coincides with $k(\GL_n/G)$.
\end{proof}

See \cite{Po5}*{Remark (1.5.7)} for a generalization of the above statement.

There also exist projective versions of the no-name method.

\begin{cor}                                                      \label{cor.noname.proj}
Let $G$ be a reductive $k$-group, let $V$
be a finite-dimensional $k$-vector space with a linear
$G$-action and let $X$ be a geometrically integral $k$-variety with an almost free
$G$-action. Then the field
$k(X\times\D{P}(V))^G$ is pure over the field $k(X)^G$. In particular, if
$k(X)^G$ is pure, so is
$k(X\times\D{P}(V))^G$.
\end{cor}

\begin{proof} 
By restricting $X$ to a
suitable open set, we may assume that $X$ is affine and that the action of
$G$ on $X$ is free.
According to the proof of Theorem \ref{thm.noname.general}, the projection $(X\times
V)\quota G\to X\quota G$ makes the
first space into a vector bundle over the second one (the action on $X\times
V$ being the diagonal
one). Since the action of $G$ commutes with that of $\D{G}_m$ on the factor
$V$, the action of $\D{G}_m$
descends to an action on $(X\times V)\quota G$, which respects the vector bundle
structure. Let $U\subset
X\quota G$ be an open set over which the vector bundle has a section, i.e.
$(X\times V)\quota G\cong U\times V$.
We now have:
\begin{align*}
k(X\times \D{P}(V))^G 
& = k(X\times V)^{\D{G}_m} = k((X\times V)\quota G)^{G\times
\D{G}_m} \\
& = k(U\times V)^{\D{G}_m} = k(U\times\D{P}(V)) = k((X\quota G)\times\D{P}(V))
\end{align*}
and this last field is clearly purely transcendental over $k(X\quota G)$.
\end{proof}

There also exists a generalized version of the no-name lemma, due to
Bogomolov \cite{Bo2}, 
and which he
has used to study fields of invariants $k(V)^G$ when $G$ is a simply connected
semisimple group.

\subsection{Combining both methods}  \label{ssec.both.method}

We now restrict attention to linear actions and prove a basic proposition
which uses both the ``no-name
lemma'' and the ``slice method''.

\begin{prop}                                            \label{prop.noname.slice}
Let $G$ be a reductive $k$-group and let
$W$ be a finite-dimensional $k$-vector space
with a linear $G$-action. Let $v\in W$ be a point such that:
\begin{enumerate}[\rm(i)]
\item
the stabilizer of $v$ is a reductive $k$-group $H$;
\item
the orbit of $v$ under $G$ is dense in $W$.
\end{enumerate}
Then if $V$ is a finite-dimensional $k$-vector space with an almost free linear
$G$-action, the action of $H$ on $V$
is almost free and the field $k(V)^H$ is pure over $k(V)^G$. In particular,
$k(V)^H$ is stably pure
if and only if $k(V)^G$ is stably pure.
\end{prop}

\begin{proof} 
Consider the diagonal action of $G$ on $W\times V$.
This is an almost free
action. Indeed, let $U\subset V$ be a $G$-invariant open set and let $x\in
U$ be such that $G_x=1$
and $G\cdot x$ is closed in $U$. Then the same properties hold for the open
set $W\times U$ and the
point $y=(0,x)\in  V\times U$. Since the action of $G$ on $V$ is almost
free, the no-name lemma
implies that $k(W\oplus V)^G$ is pure over $k(V)^G$. On the other hand, the
closed subvariety
$Y=\accol{v}\times V\subset X=W\times V$ satisfies the assumptions of
Theorem \ref{thm.slice.lemma} ---  the ``slice lemma'' ---
for the pair $(H,G)$. Indeed, $H$ is
clearly the normalizer of $Y$ and if $x$ and $y$ belong to $Y(\bar k)$ and $g\in
G(\bar k)$ satisfies $gx=y$, then
$g\in  H(\bar k)$. That the orbit of $Y$ under $G$ is dense in $W\times V$ follows
from the density of
$G\cdot x$ in $W$. Thus by the ``slice lemma'' $k(Y)^H$ is isomorphic to
$k(V\oplus W)^G$. Also, the
action of $H$ on $Y$ is isomorphic to the action of $H$ on the vector space
$V$, which is an almost
free action since the action of $G$ itself was almost free (use the same
open set $U$ and the same
point $x$ as above: since $H$ is closed in $G$, $H\cdot x$ is closed in
$G\cdot x$ hence in $U$). Thus
$k(V)^H$ is pure over $k(V)^G$.
\end{proof}

There is a projective variant of the above theorem:

\begin{prop}                                            \label{prop.noname.slice.proj}
Let $G$ be a reductive group over $k$ and let
$W$ be a finite-dimensional $k$-vector space
with a linear $G$-action. Let $v\in  W$ be a point such that:
\begin{enumerate}[\rm(i)]
\item
the stabilizer of $v$ is a reductive group $H$;
\item
the orbit of $v$ under the action of $G\times \D{G}_m$ is
dense in $W$.
\end{enumerate}
Then if $V$ is a finite-dimensional $k$-vector space with an almost free linear $G$-action, the
action of $H$ on $V$
is almost free and the field $k(\D{P}(V))^H$ is stably pure over
$k(\D{P}(V))^G$.
\end{prop}

\subsection{A criterion for retract rationality}  
                                \label{ssec.retract.rationality.criterium}

\begin{prop}                                              \label{prop.retract.rationality.criterium}
Let $G$ be a reductive
group over $k$. Assume that for any
local $k$-algebra $A$, with residue field $\kappa$, the reduction map
on \'etale cohomology sets
$H^1(A,G) \to H^1(\kappa,G)$ is onto.  Let $V$ be a $k$-vector space
with an almost free
linear $G$-action.                                     
Then $k(V)^G$ is the function field of a retract rational $k$-variety.
\end{prop}

\begin{proof}
One easily checks that if $X$ and $Y$ are two
stably $k$-birationally equivalent geometrically integral $k$-varieties,
then $X$ is retract rational over $k$ if and only if $Y$ is.
By the no-name lemma 
one is therefore reduced to
proving that if  $G \subset\GL_n$ is a closed embedding of algebraic groups,
and $X=\GL_n/G$ is the quotient variety, then
the variety $X$ is retract rational.
For $A$ and $\kappa$ as above,
we have the natural
compatible exact sequences 
of pointed sets
\[
\begin{CD}
        \GL_n(A)@>>> X(A)@>>> H^1(A,G)@>>>H^1(A,\GL_n)\\[-1mm]
      @VVV              @VVV         @VVV           @VVV\\[-1mm]
\GL_n(\kappa)@>>> X(\kappa)@>>> H^1(\kappa,G)@>>>H^1(\kappa,\GL_n).
\end{CD}
\]
Grothendieck's version of Hilbert's theorem $90$ ensures $H^1(A,\GL_n)=0$,
hence the map  $X(A) \to H^1(A,G)$ is onto.
The first arrow in each sequence is equivariant with
respect to the left action of $\GL_n(A)$, resp. $\GL_n(\kappa)$.
Two elements of $X(\kappa)$ have the same image
in $ H^1(\kappa,G)$ if and only if they are in
the same orbit of $G(\kappa)$.
If the map $H^1(A,G) \to H^1(\kappa,G)$ is onto, one concludes
that the reduction map $X(A) \to X(\kappa)$ is onto.
It remains to apply Proposition \ref{prop.retract.rational}.
\end{proof}

\section{Examples}                    \label{sec.examples}

In this section, $k$ denotes a field of characteristic zero and $\bar k$ an algebraic closure of $k$.

\subsection{Solvable Groups}  \label{ssec.solvable.gps} 

Details upon the following results may be found in the notes by Kervaire
and Vust \cite{KeVu}. 

\begin{thm}[T.~Miyata \cite{Mi}, \`E.~Vinberg \cite{V}]    \label{thm.triangular.purity}
Let the abstract
group $G$ act 
linearly
upon the
finite dimensional vector
space $V$ over $k$, and assume that there is a complete flag of $V$ which is invariant
under the action of $G$.
Then the field $k(V)^G$ is pure.
\end{thm}

See also \cite{Tr}.

\begin{proof} 
The proof is easily reduced to the following nice
lemma, whose proof only
uses the euclidean division algorithm.
\end{proof}

\begin{lem}                                    \label{lem.invariant.descent}
Let $K$ be a field, and let the abstract
group $G$ act upon the
polynomial ring in one variable $K\croch{t}$. If $K\subset  K\croch{t}$ is globally
fixed, then there exists a
$G$-invariant polynomial $p$ such that $K(t)^G=K^G(p)$.
\end{lem}

\begin{prop}                                                             \label{prop.torus.purity}
Assume $k=\bar k$.
If $V$ is a finite dimensional
vector  space over $k$ and $G \subset
 \GL(V)$ is an
(abstract) abelian group consisting of semisimple elements, then
$k(V)^G$ is pure.
\end{prop}

For $G$ a finite group this is just Fischer's well known theorem \cite{Fi}.

\begin{proof} 
Indeed, all elements of $G$ may be simultaneously
diagonalized, so that $G$
actually acts through a maximal torus of $\GL(V)$.
\end{proof}

\begin{prop}                                  \label{prop.triangular.purity}
Assume $k=\bar k$. 
If $G$ is a connected solvable
group over $k$
and if $G$ acts linearly
on the finite dimensional
vector  space $V$ over $k$, then $k(V)^G$ is pure.
\end{prop}

\begin{proof} 
This is a consequence of the above theorem and of the
Lie-Kolchin theorem.
\end{proof}

If the (abstract group) $G$ acts linearly on the vector space $V$, this
action induces an action on
the projective space $\D{P}(V)$.

\begin{prop}                                      \label{prop.invariant.descent}
The field $k(V)^G$ is pure over
$k(\D{P}(V))^G$.
\end{prop}

\begin{proof} 
This follows from Lemma \ref{lem.invariant.descent}.
\end{proof}

For more general actions, we have the following

\begin{thm}[Rosenlicht \cites{R1,R3}]                      \label{thm.rosenlicht}
Let
$G$ be a connected solvable algebraic group over an algebraically closed field $k$.
If $X$ is an integral variety
with a $G$-action, and
$f\colon X\to Y$ is a geometric quotient 
for this action, then $f$ admits a
section over a non-empty open set
of $Y$.
\end{thm}

This theorem has the following consequences:
\begin{enumerate}[(i)]
\item
If $X$ is a rational variety, then $Y$ is a retract rational variety.
\item
If the action of $G$ on $X$ is free, i.e.\ if $X$ is a $G$-torsor over
$Y$, then $X$ is birational
to $G\times Y$; if moreover $X$ is rational, then $Y$ is stably rational.
\end{enumerate}

\subsection{Special Groups}       \label{ssec.special.gps}

Recall that an algebraic group $G$ over a field $k$ 
is called \emph{special}
if any principal
homogeneous space (torsor) under
$G$ over a $k$-variety  is locally trivial for the Zariski topology.
This implies $H^1(K,G)=0$  that for any field $K$ containing $k$.

The following facts were
proved by Serre \cite{Ch}.
Any
special group $G$ is linear and connected. If $G\subset \GL_n$ is some
closed embedding, $G$ is
special if and only if the fibration $\GL_n\to\GL_n/G$ is locally trivial
for the Zariski topology.

Special groups include 
\begin{enumerate}[\rm(i)]
\item
the additive group $\D{G}_a$,
\item
 the multiplicative
group $\D{G}_m$,
\item
 more generally, split connected solvable groups,
 \end{enumerate}
as well as some connected linear algebraic groups, e.g., 
\begin{enumerate}[\rm(i)]
\addtocounter{enumi}{4}
\item
the
general linear groups
$\GL_n$, more generally, $\GL(A)$ for $A$ a central simple algebra over $k$,
\item
the
special linear groups
$\SL_n$, 
\item
the symplectic groups $\Sp_{2n}$, 
\item
the split spinor groups $\Spin_n$ for
$n\leq 6$.
\end{enumerate}
This last case
follows from the existence of the exceptional isomorphisms in the
classification of Lie groups: $\Spin_3\cong\SL_2$, $\Spin_4\cong\SL_2\times\SL_2$, $\Spin_5\cong\Sp_4$, 
$\Spin_6\cong\SL_4$. The
special orthogonal groups $\SO_n$ for $n\geq 3$ and the spinor groups for
$n\geq 7$ are not special.

Over an algebraically closed field, any
connected linear algebraic group is a rational variety. 
Over an arbitrary field $k$, if a  semisimple group  is special, it is a product  of copies of Weil restriction of scalars of groups of type $\SL_n$ and $\Sp_{2n}$
(Grothendieck \cite{Ch} over an algebraically closed field, one checks this extends to arbitrary $k$).
Thus the underlying variety of a special semisimple $k$-group is a $k$-rational variety. This statement is not true for special $k$-tori. 

\begin{prop}                                                \label{prop.special.stable.purity}
Let $G$ be a special group over $k$
whose underlying variety is $k$-rational, and let
$X$ be a factorial affine $k$-variety
with an almost free $G$-action. Then $k(X)$ is purely transcendental over
$k(X)^G$. In particular,
if $X$ is a (stably) rational variety, then $k(X)^G$ is stably pure.
\end{prop}

\begin{proof} 
Since the action is almost free, there
exist  a non-empty open set $U\subset X$ and a morphism $U\to Y$ which makes $U$
into a $G$-torsor over $Y$ (Corollary \ref{cor.vdbergh.reductif}). In
particular, the function field $k(Y)$ coincides with $k(X)^G$. Because
$G$ is special, any
$G$-torsor is locally trivial for the Zariski topology and $U$, hence also
$X$, is birational to the
product $Y\times G$. Thus $k(X)$ is pure over $k(Y)=k(X)^G$.
\end{proof}

\begin{cor}                                                                \label{cor.special.stable.purity}
Let $G$ be a semisimple special group over $k$, 
and let $V$ be a finite dimensional $k$-vector space equipped with an almost free linear
$G$-action. 
Then $k(V)^G$ is stably pure.
\end{cor}

\begin{proof} 
Indeed, as mentioned above, $G$ is a $k$-rational variety.
\end{proof} 

\begin{prop}                                                \label{prop.special.stable.equivalence}
Let $G$ be a linear 
algebraic group over
$k$. Let $G_1,G_2$ be two special groups 
which contain $G$ as a closed
subgroup.
Assume that  the underlying varieties of $G_1$ and $G_2$ are rational over $k$.
Then $G_1/G$ and $G_2/G$ are stably $k$-birationally
equivalent.
\end{prop}

\begin{proof}
Let $X=(G_1\times_kG_2)/G$ be the quotient of 
$G_1 \times_kG_2$ under  the diagonal
embedding. The projection map $X \to G_1/G$
makes $X$ into a $G_2$-torsor over $G_1/G$.
Since $G_2$ is special, $X$ is $k$-birational to
the product $G_2\times (G_1/G)$. 
Thus $X$ is stably $k$-birational to $G_1/G$. The same argument
applies to the projection $X \to G_2/G$. Hence
$G_1/G$ is stably $k$-birationally equivalent to
$G_2/G$.
\end{proof}

\subsection{The Orthogonal Group}                 
                                                      \label{ssec.orth.group}
                                                 
\begin{prop}                                                                   \label{prop.orthogonal.purity}
Let $q$ be a nondegenerate quadratic form over $k$ of rank
$n\geq 2$. Let $G=\Ort(q)$ be
the orthogonal group of $q$. 
Let $V$ be a finite-dimensional $k$-vector space with an almost free 
linear $G$-action. Then the field $k(V)^G$ is stably pure over $k$.
\end{prop}

\begin{proof}[1st Proof (by reduction to $\GL_n$)]
Apply Proposition \ref{prop.noname.slice}
to the group $G=\GL_n$
acting upon the vector space $V$ of symmetric matrices $M$ through $g\cdot
M={}^tgMg$, and
take $v$ in that proposition to be the matrix $S$ of $q$. The orbit of $v$ is
the open set of
non-singular symmetric matrices, hence it is dense. The stabilizer of $v$
is
$\Ort(q)$. Then Proposition \ref{prop.noname.slice} reduces the case of $\Ort(q)$ to that of
$\GL_n$ which is a
``special'' group. Then we are done by Proposition \ref{prop.special.stable.purity}.
\end{proof}

\begin{proof}[2nd Proof]
Let $Q_n\subset \C{M}_n$ be the vector space
consisting of symmetric
matrices. Consider the map $\C{M}_n\to Q_n$ which sends the matrix $A$ to
${}^tASA$. This map
induces a map $f$ of affine varieties $X=\GL_n\to Y=Q_n^{\text{ns}}$, where  
$Q_n^{\text{ns}}$ denotes
the open set of
non-singular symmetric  $n\times n$ matrices. Let $G=\Ort(q)$ act upon $\C{M}_n$, hence
also $\GL_n$, by left
translation, i.e.\ by the action $(g,A)\mapsto gA$. One easily checks that
$f$ makes $X$ into a
$G$-torsor over $Y$ for this action, hence $k(X)^G=k(Y)$ (that the
map $f$ is onto over $\bar k$ reflects the
fact that any non-degenerate quadratic form of rank $n$ is equivalent to $q$ over $\bar k$). Since
$Y=Q_n^{\text{ns}}$ is an open set in a vector space, $k(Y)$ is purely transcendental. On
the other hand, the
action of $G$ on $\C{M}_n$ is linear, and it is an almost free action. The
theorem now follows from
the no-name lemma (Corollary \ref{cor.noname.lemma}).
\end{proof}

\begin{rem}                                         \label{rem.orth.retract.rationality}
Using Proposition \ref{prop.retract.rational}, 
the retract rationality of
$k(V)^G$ can already be seen
from the following fact: \emph{For any local $k$-algebra $A$ with residue class
field $\kappa$, the natural map
$H^1(A,\Ort(q))\to H^1(\kappa,\Ort(q))$ is surjective}. Indeed, $H^1(\kappa,\Ort(q))$
classifies isomorphism classes of
non-degenerate rank $n$ quadratic forms over $\kappa$. Any such class $\alpha$ contains
a diagonal form $a_1x_1^2+\dots+a_nx_n^2$. Lifting the $a_i\in
\kappa^*$ to elements of
$A^*$ produces a non-degenerate quadratic form, hence an element in
$H^1(A,\Ort(q))$ whose image under
the restriction map is $\alpha$.
\end{rem}

\subsection{The Special  Orthogonal 
                        Group}         \label{ssec.special.orth.group}

Let $q$ be a nondegenerate quadratic form over $k$ of rank
$n\geq 2$, let $\disc(q)\in k^*/k^{*2}$ be its discriminant. Let $G=\SO(q)$ be the special orthogonal group of $q$. 
If $A$ is a local $k$-algebra
with $2\in  A^*$, then $H^1(A,\SO(q))$ classifies non-degenerate quadratic
forms over $A$ with the same
discriminant as $q$ in $A^*/A^{*2}$. Now, if $A$ is a local ring with residue class
field $\kappa$ ($\car\kappa\neq 2$)
and $\sum a_ix_i^2$ is a non-degenerate quadratic form over $\kappa$ with
discriminant $\disc(q)$, it can clearly be lifted to a non-degenerate
quadratic form over $A$ with
the same property. 

Just as in Remark \ref{rem.orth.retract.rationality}, the lifting property just proven
implies that for an almost 
free linear $G$-action of
$G=\SO(q)$ on a finite dimensional vector space $V$ over $k$ the field $k(V)^G$ is retract rational. 
One can actually do better.

\begin{prop}                                     \label{prop.rotation.group.purity}
Let $q$ be a nondegenerate 
quadratic form over $k$ of rank
$n\geq 2$. For
any almost free finite dimensional linear representation $V$ 
of the special orthogonal group $\SO(q)$ over $k$, the field
$k(V)^{\SO(q)}$ is stably rational.
\end{prop}

\begin{proof}
 Let 
$X=V^{n-1}$, let $Y=k^N$, $N=n(n-1)/2$, and let $f\colon X\to Y$ be the map
which sends
$(e_1,\dots,e_{n-1})$ to the point $\accol{\scal{e_i,e_j}}_{1\leq i\leq j\leq
n-1}$, where $\scal{x,y}$
denotes the ``scalar product'' of two vectors $x$ and $y$ with respect to the bilinear form assoiated to $q$. Applying 
Proposition
\ref{prop.generic.quotient} to the open set $U$ of $X$ consisting of $(e_1,\dots,e_{n-1})$ with
$e_1\wedge\dots\wedge e_{n-1}\neq 0$, one obtains
$k(Y)=k(X)^G$. Moreover, the action of $G$ on $X$ is almost
free. Since $Y$ is
rational, this produces an almost-free finite dimensional linear representation $X$ over $k$ such that $k(X)^G$ is pure over $k$. We conclude by the no-name lemma (Corollary \ref{cor.noname.lemma}).
\end{proof}

\begin{rem}                                                      \label{rem.bogo.saltman.proofs}
To prove this result, Bogomolov \cite{Bo2}
uses a generalized version
of the no-name lemma
 together with
the slice method; see also Saltman's approach \cite{Sa8}.
\end{rem}

\subsection{The Spinor Group}     \label{ssec.spin.group}
                                                      
Let $\Spin(q)$ denote the spinor group attached to
a nondegenerate quadratic form $q$ of rank $n$.
If $k=\bar k$, we write simply $\Spin(q)=\Spin_n$.

Assume $k=\bar k$. As we saw in \S\ref{ssec.special.gps} for $n\leq 6$, the group $\Spin_n$ is special, hence if
$V$ is an almost free finite dimensional linear representation of $\Spin_n$,
the field $k(V)^{\Spin_n}$ is stably pure.
In the literature, there is  a proof \cite{Bo2} of the stable purity of $k(V)^{\Spin_n}$ for arbitrary $n$, 
but that proof is incomplete, it builds upon an incorrectly proved lemma \cite{Bo2}*{Lemma 3.1}.
A proof of the stable rationality of $\Spin_n$ for $n=7$ and $n=10$ is given by 
V.\`E.~Kordonski\u{\i} \cite{Kor1}.

D.~Shapiro suggested that the classification of
forms of low degree given by A.~Pfister in 1966 would immediately
yield the retract rationality of $k(V)^{\Spin_n}$  for $n\leq 12$.
Here we shall only discuss the case $n=12$.

\begin{thm}                                           \label{thm.spin.retract.rationality}
Let $G=\Spin(q_0)$ be
the spinor group of the nondegenerate hyperbolic quadratic form $q_0$ of rank $12$.
Let $V$ be a finite-dimensional $k$-vector space with an almost free 
linear 
$G$-action. 
Then the field $k(V)^G$ is retract rational over $k$.
\end{thm}

\begin{proof}
Let us write $\Spin=\Spin(q_0)$ and $\SO=\SO(q_0)$.
Let
\[
1\to\mu_2\to\Spin\to\SO\to 1
\]
be the natural sequence. For any $k$-algebra $A$,
there is an an exact sequence of pointed sets in \'etale cohomology:
\[
H^1(A,\mu_2) \to H^1(A,\Spin) \to H^1(A,\SO)\xrightarrow{\partial}H^2(A,\mu_2).
\]
There is an action of the group $H^1(A,\mu_2)$ on
$H^1(A,\Spin)$, and this action is transitive on elements of
$H^1(A,\Spin)$ with the same image in $H^1(A,\SO)$ 
\citelist{\cite{Ch}\cite{Se2}*{I. 5.7, Prop.~42 
p.~52}}.
The set $H^1(A,\SO)$ classifies the non-degenerate quadratic forms
over $A$ of dimension $n=12$ with discriminant $1\in A^*/A^{*2}$.
The map  
\[
H^1(A,\SO)\xrightarrow{\partial}H^2(A,\mu_2)
\]
associates to the
class of a such a form its Clifford invariant \citelist{\cite{Se2}*{pp.~147\textendash 148}
\cite{KMRT}*{p.~437}}.

Let $A$ be a local $k$-algebra, let $\kappa$ be its residue field.
Let $\xi\in H^1(\kappa,\Spin)$. Its image in $H^1(\kappa,\SO)$
is the isomorphy class of a quadratic form $q$ of rank $12$, trivial
discriminant, and trivial Clifford invariant.
According to a theorem of Pfister \cite{Pf}*{Satz 14} (see also \cite{Kah})
there
exist elements $a, b, c, d, e, f$ in $k^*$ such that
such a form may be written as
$q= a \scal{1,b}\otimes\scal{-c,-d,cd, e,f, -ef}$.
Since the reduction map $A^*/A^{*2} \to \kappa^*/\kappa^{*2}$ is onto,
one may lift this quadratic form to a class of the same
shape over $A$, and the Clifford invariant of such a class
is trivial, as may be checked by a direct computation.                
Indeed, such a class clearly lies in $I^3A$.                                 
One thus finds an element $\eta_A \in H^1(A,\Spin)$
whose image  $\eta_{\kappa} \in H^1(\kappa,\Spin)$ has the same image
as $\xi$ in  $H^1(\kappa,\Spin)$. There
thus exists $\rho_{\kappa} \in H^1(\kappa,\mu_2)$
such that $\rho_{\kappa}.\eta_{\kappa}=\xi$.
Let $\rho \in H^1(A,\mu_2)$ be a lift of this element.
Then $\rho\cdot\eta \in H^1(A,\Spin)$ reduces to 
$\xi \in H^1(\kappa,\Spin)$.                                                          
\end{proof}

\subsection{The Projective Linear Group}     
                                                 \label{ssec.projective.group}
                                                 
Let $V$ be a finite dimensional vector space over $k$,
$n\geq 2$ an integer. Suppose the projective linear group $\PGL_n$ acts linearly and almost freely on $V$.
The question whether the field of invariants $k(V)^{\PGL_n}$ is rational (pure), or at least stably rational,
has come up in a variety of contexts. It is still open.
 
The simplest linear representation of the above type is
 given by $V_n=\C{M}_n(k)\oplus\C{M}_n(k)$ with 
 $\PGL_n$ acting by simultaneous conjugation.
 
The field of invariants $C_n: = k(V_n)^{\PGL_n}$ coincides with
the field of fractions of the centre of the generic division ring on two $n\times n$ generic matrices
(Procesi \cite{Pro1}*{pp.~240\textendash 241};  M.~Artin \cite{Ar}).  

At least for $k$ the field of complex numbers, the field $C_n$ may also be viewed as 
\begin{enumerate}[---]
\item
The function field 
of the moduli space $M(n;0,n)$ of stable rank $n$ vector bundles on $\D{P}_2$ with Chern numbers $c_1=0$, $c_2=n$ (see Le Bruyn \cites{LeBr1,LeBr2}).
\item
The function field of the generic Jacobian variety of plane curves of degree $n$
(Van den Bergh \cite{VdB}).
\end{enumerate}
Katsylo and Schofield proved the following reduction result (see also \cite{Sa11}):

\begin{prop}[\cites{Kat6, Scho}]                                                \label{rem.stably.equivalence.coprime}
Suppose $n=rs$ with $\gcd(r,s)=1$, 
then $C_n$ is stably equivalent to 
the fraction field of $C_r\otimes C_s$.
\end{prop}

Saltman proved retract rationality of $C_n$ for $n$ prime \cite{Sa2}*{Corollary 5.3}. 
For an alternate proof, see also \cite{CS2}*{Corollary 9.13}.

\begin{prop}                                     \label{prop.center.generic.dva.retract.rationality}
Let $n$ be squarefree or twice a squarefree number. Then for any almost free
finite dimensional linear representation 
$V$ of the projective linear group group $\PGL_n$, the field
$k(V)^{\PGL_n}$ is  retract rational.
\end{prop}

The rationality for $n=2$ was proved by Sylvester \cite{Sy}, and by Procesi \cite{Pro1}.
 Formanek solved the rationality problem for $n=3$ \cite{Fo1} and $n=4$ \cite{Fo2}. 
 For $n$ arbitrary, Formanek \cite{Fo1} gave a very useful description
of $C_n$ as a field of multiplicative invariants under an action of the symmetric group.
For $k=\bar k$, Bessenrodt and Le Bruyn \cite{BL} proved the stable rationality for $n=5,7$.
A simpler proof was later given by Beneish \cite{Be1}.
The question remains open for all other prime powers $n$.
The quoted results combine to :

\begin{prop}                                     \label{prop.center.generic.divalg.purity}
Let $k=\bar k$. Let $n=2,3,4,5,7$, or more generally 
let $n$ be any divisor of $420$.
Then for any almost free linear representation 
$V$ of the projective linear group group $\PGL_n$, the field
$k(V)^{\PGL_n}$ is 
stably rational.
\end{prop}

For further rationality results on $k(V)^{\PGL_n}$, with connections with the rationality of some moduli varieties,
see \cites{Bo1,BoKat,BoPT,Kat1,Kat2,Kat3,Kat4,Kat5,Kat6,Kat7,Kat8,Kat9,Kat10,Kat11,ShB1,ShB2,ShB3,ShB4}.
The case of some exceptional groups is considered in \cites{Ilt,IltSh}.

\subsection{Some Finite Groups}  \label{ssec.finite.gps}

\begin{lem}                                 \label{prop.finite.group.free}
Let $X$ be a 
geometrically  integral $k$-variety, and
$G$ a finite group acting on $X$. If the action of $G$ is faithful, then the action is almost free.
\end{lem}

\begin{proof} 
For any $g \in G$, $g\neq e$,
the $k$-variety $X^g \subset X$ of fixed points has codimension
at least $1$. Let $U$ be the complement
of the union of the $X^g$,  $g\neq e$. Then $U$ is a $G$-stable nonempty open
set and the action of $G$ on $U$ is free.
\end{proof} 

\begin{prop}                                     \label{prop.permutation.group.purity}
For
any faithful finite dimensional linear representation $V$ of $\F{S}_n$ over $k$, the field
$k(V)^{\F{S}_n}/k$ is stably pure over $k$.
\end{prop}

\begin{proof} 
This is the most classical case. For the
natural linear representation of $\F{S}_n$ in $k^n$ by permutations on the
canonical basis, the theorem on symmetric polynomials says that
\[
k(t_1,\ldots,t_n)^{\F{s}_n}=k(s_1,\ldots,s_n)
\]
is a pure extension of $k$. Here we denote by
$s_1=t_1+\dots+t_n,\ldots,s_n=t_1\dots t_n$ the
fundamental symmetric polynomials.
The proposition now follows from the no-name lemma.
\end{proof}

But the analogous question for $\F{A}_n$, $n \geq 6$  remains open!

\begin{prop}                                    \label{prop.finite.subgroup.gl.three}
Assume $k=\bar k$.
Let $G$ be a finite subgroup of
$\GL(3,k)$. Then for
any faithful finite dimensional linear representation $V$ of $G$ over $k$, the field  $k(V)^G/k$ is stably pure
over $k$.
\end{prop}

\begin{proof}
Let $V_0$ be a faithful representation, $\dim V\leq 3$.
By Proposition \ref{prop.invariant.descent} we know that the field
$k(V_0)^G$ is pure over the field $k(\D{P}(V_0))^G$. Since this last field is
unirational of transcendence
degree at most $2$, it is purely transcendental by Castelnuovo's theorem. The
result now follows from
the no-name lemma.
\end{proof}

\begin{rem}                                                        \label{rem.finite.subgroup.gl.two}
Let  $k$ be a field and $G\subset\GL(2,k)$ be any finite subgroup.
Then for any faithful finite dimensional linear representation $V$ of $G$ over $k$, the field $k(V)^G/k$ is stably pure. This result applies to the dihedral group $G=D_n$ of order $2n$ and $k=\bar k$ or $k=\D{R}$.
For the proof we consider a faithful representation $V_0$, $\dim V_0\leq 2$. By Proposition \ref{prop.invariant.descent} and L\"uroth's theorem, $k(V_0)^G/k$ is pure. Hence $k(V)^G/k$ is stably pure by the no-name lemma.
\end{rem}

In some cases, one may use the no-name lemma to prove the purity of fields
of invariants.

\begin{thm}[Maeda \cite{Ma}]                               \label{thm.fifth.alternating.group.purity}
Let $V$ be the natural
$5$-dimensional
representation of the alternating group $\F{A}_5$.
Then $k(V)^{\F{A}_5}$
is pure over $k$.
\end{thm}

We here give a proof in the case $k=\bar k$.

\begin{proof}[Proof (for $k=\bar k$)]
There exists a faithful representation $W$ of
dimension 3 (as the group of
automorphisms of the icosahedron). Also, there is a decomposition
$V=V_0\oplus V_1$ where $V_0$ is
the trivial representation. Now $V/\F{A}_5$ is pure over
$\D{P}(V_1)/\F{A}_5$ and of transcendence degree
$2$. On the other hand, the action of $\F{A}_5$ on $\D{P}(V_1)$ and $\D{P}(W)$ is
almost free (indeed,
$\F{A}_5$ is simple). 
We may thus apply the projective no-name lemma, i.e.\ Corollary \ref{cor.noname.proj}. This
implies that
$(\D{P}(V_1)\times \D{P}(W))/\F{A}_5$ is pure over $\D{P}(V_1)/\F{A}_5$ of
transcendence degree $2$. Thus
$\D{P}(V_1)/\F{A}_5$ and $(\D{P}(V_1)\times \D{P}(W))/\F{A}_5$ are birational. A
second application of the
projective no-name lemma implies that $(\D{P}(V_1)\times \D{P}(W))/\F{A}_5$ is
pure over
$\D{P}(W)/\F{A}_5$, which is unirational of dimension $2$, hence pure by Castelnuovo's theorem
since $\car k=0$.
\end{proof}

We refer to \cite{Ma} for the proof over $\D{Q}$. In \cite{Pl}, B.~Plans shows that for odd $n\geq 3$, the field 
$\D{Q}(X_1,\dots,X_n)^{\F{A}_n}$ is pure over $\D{Q}(X_1,\dots,X_{n-1})^{\F{A}_{n-1}}$, thus giving another proof of Maeda's theorem.

\begin{rem}                                           \label{rem.fifth.alternating.group.purity}
One may also prove that 
$k(\D{P}(V))^{\F{A}_5}$ is pure, see \cite{HaKan3}*{Lemma 5}.
\end{rem}

\begin{rem}                                           \label{rem.icosahedral.group.purity}
The binary icosahedral group 
is a subgroup of
$\SL_2(\D{C})$. 
It is a double cover of $\F{A}_5$, actually it is its ``representation group'' in the sense of Schur.
It has a unique  $4$-dimensional faithful linear action.
The analysis of the quotient 
 requires a much more delicate analysis \cite{KoPr1}. 
See also \cites{KoPr2,KoPr3} for some other finite subgroups 
of $\GL_{4}$.
\end{rem}

For further results on the rationality of the field of invariants for a linear action
of other small finite groups, see \cites{AHK,ChHKan,ChKan,Sa7,Prok}.

\section{The unramified Brauer group}
                                                                  \label{sec.def.bru}

In this section, we shall assume some knowledge of group
cohomology, as well as some
knowledge of local fields, all of which may be found in Serre's book \cite{Se1}.
For the sake of
simplicity, all fields ($k,K,L,\kappa,\dots$) will be taken of characteristic zero. Let $K$ be a
field, let $\overline K$  be an
algebraic closure of $K$ and let $\F{g}=\Gal(\overline K/K)$.

\begin{Def}                                                    \label{def.brauer.group}
The second (profinite) cohomology group
$H^2(\F{g},\overline K^*)$ of
$\F{g}$ with values in the multiplicative group $\overline K^*$ of
$\overline K$ is called the Brauer
group of $K$ and is denoted $\Br K$.
\end{Def}

Given any field inclusion $K\subset  L$, there is a natural map $\Br K\to
\Br L$.

\begin{Def}                                                    \label{def.residue.maps}
When $K$ is the field of fractions
of a discrete valuation ring $A$ with residue field $\kappa$ (of
characteristic zero), there is a
basic homomorphism
\[
\partial_A\colon\Br K \to \Char(\kappa),
\]
where $\Char(\kappa)$ denotes the group $\Hom_{\text{cont}}(\Gal(\bar
\kappa/\kappa),\D{Q}/\D{Z})$ of continuous
characters of $\Gal(\bar\kappa/\kappa)$ with values in the discrete group $\D{Q}/\D{Z}$.
\end{Def}

Let us first assume that $A$ is complete. Since the characteristic of
$\kappa$ is assumed to be
zero, $A$, resp.\ $K$, may be then identified with the ring of power series
$\kappa[[t]]$, resp.\ with
its fraction field $\kappa((t))$. Let $K_{\text{nr}}$ be the maximal
unramified extension of $K$. Under
the previous identifications, $K_{\text{nr}}$ coincides with $\bar k((t))_{\text{alg}}$
(the algebraic closure of $\kappa((t))$ in $\bar k((t))$),
$\Gal(K_{\text{nr}}/K) = \Gal(\bar\kappa/\kappa)$, and $\overline K$ is the union
of all fields
$\bar\kappa((t^{1/n}))_{\text{alg}}$ (``Puiseux's theorem''). The Galois group
$\Gal(\overline K/K_{\text{nr}})$ can thus be identified with the profinite group  $\varprojlim\mu_n$,
which in a non
canonical manner is isomorphic to the group $\widehat{\D{Z}}=\varprojlim \D{Z}/n$.
The cohomological
dimension of such a group is $1$, hence $H^2(\Gal(\overline K/K_{\text{nr}}),\overline K^*)=0$. 
Since also
$H^1(\Gal(\overline K/K_{\text{nr}}),\overline K^*)=0$ according to Hilbert's
Theorem $90$, the
restriction-inflation sequence
gives rise to a natural isomorphism:
\begin{equation}                                                               \label{eq.brauer.isomorphism.nr}
H^2(\Gal(K_{\text{nr}}/K),K_{\text{nr}}^*) 
\xrightarrow{\cong}H^2(\Gal(\overline K/K),\overline K^*)=\Br K.
\end{equation}
The
valuation $v\colon K^*\to\D{Z}$ naturally extends to a valuation $v\colon K_{\text{nr}}^*
\to\D{Z}$ which respects the
action of $\Gal(K_{\text{nr}}/K)$ ($\D{Z}$ being taken with trivial action). It
therefore induces a
homomorphism:
\begin{equation}                                                          \label{eq.brauer.morphism.residu}
H^2(\Gal(K_{\text{nr}}/K),K_{\text{nr}}^*) \to  H^2(\Gal(K_{\text{nr}}/K),\D{Z}) =
H^2(\Gal(\bar\kappa/\kappa),\D{Z}).
\end{equation}
Now the cohomology of the exact sequence $0 \to\D{Z} \to
\D{Q} \to\D{Q}/\D{Z} \to 0$ identifies the last group with
$H^1(\Gal(\bar\kappa/\kappa),\D{Q}/\D{Z})
=\Hom_{\text{cont}}(\Gal(\bar\kappa/\kappa),\D{Q}/\D{Z})$.
Combining the isomorphism \eqref{eq.brauer.isomorphism.nr}
with the map \eqref{eq.brauer.morphism.residu} defines the
map $\partial_A$ in the case when $A$ is complete. 

In the general case,
$\partial_A$ will be defined as the composite map
\[
\Br K \to\Br\widehat K\xrightarrow{\partial_{\hat A}}\Char(\kappa)
\]
where $\widehat A$, resp.~$\widehat K$, are the completions of $A$, resp.~$K$ --- the residue class
field $\kappa$ being the same as that of $A$.

We are now in a position to define our basic invariant, first used efficiently by Saltman \cite{Sa3}.

\begin{Def}                                                                               \label{def.bru.group}
Let $k$ be a field of characteristic zero. Let $K/k$ be a function field, i.e.\ assume
that the field
$K$ is finitely generated, as a field, over the field $k$. The subgroup
$\bigcap_A \ker\partial_A$ of
$\Br K$, where $A$ runs through all discrete rank one valuation rings
$A\subset  K$ such that the
field of fractions of $A$ is $K$ and $k$ is included in $A$, is called the
\emph{unramified Brauer
group} of $K$ (with respect to $k$) and it is denoted $\Bru(K/k)$, or $\Bru K$  if
there is no ambiguity on $k$.
\end{Def}

\begin{rem*}
In the case where $k$ is algebraically closed, $k^*$ is infinitely
divisible, hence for any 
 discrete valuation ring $A\subset K$ with fraction field $K$,
 we have $k \subset A$.
Hence, in this case, we just write $\Bru K$.
\end{rem*}

\begin{lem}                                                                              \label{lem.brk.bruK}
The natural map $\Br k\to\Br K$ sends
$\Br k$ to $\Bru(K/k)$.
\end{lem}

\begin{proof} 
It is enough to prove the statement when $K$ is the
field of fractions of a
complete discrete valuation ring $A$ containing $k$, i.e.\ $A=\kappa[[t]]$
and $K=\kappa((t))$, with
$k\subset \kappa$. The composite map
\[
\bar k^* \to\bar\kappa((t))_{\text{alg}}^* = K_{\text{nr}}^* \to\D{Z},
\]
where the last map is given by the valuation, is zero. Hence the composite
map
\[
H^2(\Gal(\bar k/k),\bar k^*) \to  H^2(\Gal(K_{\text{nr}}/K),K_{\text{nr}}^*)
\to  H^2(\Gal(\bar\kappa/\kappa),\D{Z})
=\Char(\kappa)
\] 
is zero.
\end{proof}

\begin{lem}                                                       \label{lem.functoriality.bru}
Let $K\subset  L$ be function fields over
the field $k$. The natural
map $\Br K\to\Br L$ induces a map $\Bru(K/k)\to\Bru(L/k)$.
\end{lem}

\begin{proof} 
Let $\alpha$ be an element of $\Bru(K/k)$, and let
$\alpha_L$ be
its image in $\Br L$. Let $B\subset  L$ be a discrete valuation ring of
rank one with field of
fractions $L$ and with $k\subset  L$. Let $A=B\cap K$ be the trace of $A$
on $K$.
We have $k \subset A$.
If $A=K$, then
$\partial_B(\alpha_L)=0$ according to the previous 
lemma.
Otherwise $A$ is
a discrete valuation ring
of rank one with field of fractions $K$, 
and there is a natural inclusion
of the residue class field
$\kappa_A$ of $A$ into the residue class field $\kappa_B$ of $B$. Let $\pi$
be a uniformizing
parameter of $A$, and let $e=v_B(\pi)>0$ be the valuation of $\pi$ in $B$.
The result now follows
from the general fact that in such a situation, there is a commutative
diagram:
\begin{equation}                                            \label{eq.functoriality.residue}
\begin{CD}
        \Br L@>\partial_B>>         \Char(\kappa_B)\\[-1mm]
                                       @A\Res_{K/L}AA      @AA e\cdot\Res_{\kappa_A/\kappa_B}A\\[-1mm]
       \Br K                        @>\partial_A>>        \Char(\kappa_A)
\end{CD}
\end{equation}
where the L.H.S.\ vertical map is the
natural map, and the R.H.S.\ vertical map is $e$ times the map induced on
character groups by the
inclusion $\kappa_A\subset\kappa_B$. To prove this last fact, one
simply goes over to
completions, in which case the inclusion $K\subset  L$ 
(resp.\ $K_{\text{nr}}\subset L_{\text{nr}}$)
reads $\kappa_A((t)) \subset   \kappa_B((u))$ (resp.\ $\bar\kappa_A((t))_{\text{alg}}
\subset
\bar\kappa_B((u))_{\text{alg}}$) with $u=\rho t^e$ for some $\rho\in\kappa_B^*$.
\end{proof}

\begin{lem}                                              \label{lem.bru.purity.invariance} 
If $K=k(t)$ is the
rational field in one variable over 
the field 
$k$, the natural map $\Br (k)\to\Bru(K/k)$ is an
isomorphism.
\end{lem}

\begin{proof}
Let $\bar k$ be an algebraic closure of $k$.
Let $G=\Gal(\bar k/k)$. Since $\bar k\croch{t}$ is a unique
factorization domain and its units
are reduced to $\bar k^*$, there is an obvious exact sequence of
$G$-modules:
\begin{equation}                                       \label{eq.projective.line}
1\to\bar k^*\to\bar k(t)^*\xrightarrow{\accol{v_P}}\oplus_P \D{Z}_P\to 0.
\end{equation}
Here $P$ runs through all monic irreducible polynomials in $k\croch{t}$.
If $P$ is such a polynomial, and $P(t)=(t-\alpha_1)\dots(t-\alpha_n)$ is
its decomposition over
$\bar k$, one lets $\D{Z}_P=\D{Z}\alpha_1\oplus\dots\oplus\D{Z}\alpha_n$ be the
permutation $G$-lattice
(free $\D{Z}$-module) with $\D{Z}$-basis the roots $\accol{\alpha_1,\dots,\alpha_n}$
of $P$ in $\bar k$, 
the action of $G$
being given by the permutation action on these roots. The map $\accol{v_P}$
sends any function to its
divisor, i.e.\ to the set of zeros and poles with multiplicities.

Sequence \eqref{eq.projective.line}
is a split $G$-sequence: indeed, the map which associates to $f\in\bar
k(t)^*$ the value at zero of
$f/t^{v(f)}$, where $v$ is the valuation associated to the polynomial $t$,
defines a $G$-retraction
of the embedding $k^*\to\bar k(t)^*$. Thus the $G$-cohomology of sequence
\eqref{eq.projective.line} gives rise to a
(split) short exact sequence
\begin{equation}                                        \label{eq.htwo.projective.line}
0\to  H^2(G,\bar k^*)\to  H^2(G,\bar k(t)^*)\xrightarrow{\accol{v_P}}\oplus_P
H^2(G,\D{Z}_P)\to 0.
\end{equation}
If
$k(\alpha)=k\croch{t}/P$, the group $H^2(G,\D{Z}_P)$ may be identified with
\[
H^2(\Gal(\bar k/k(\alpha),\D{Z}) =
H^1(\Gal(\bar k/k(\alpha),\D{Q}/\D{Z}) = \Char(k(\alpha)).
\]
Tsen's theorem asserts that the field $\bar k(t)$ is a ${C_1}$-field,
and, as such, its Brauer group is trivial. Using Hilbert's theorem 90 and
the restriction-inflation
sequence we may conclude that the inflation map is an isomorphism 
\[
H^2(G,\bar k(t)^*)\xrightarrow{\cong}\Br (k(t)).
\] 
Let now
$P$ be a monic irreducible polynomial in $k\croch{t}$. Such a polynomial defines
a valuation on $k(t)$. Let
$A$ be the associated discrete valuation ring, whose residue class field
may be identified with the
field $k(\alpha)=k\croch{t}/P$. By comparing with the above definition of
$\partial_A$, one then easily
checks that the following diagram is commutative:
\[
\begin{CD}
        H^2(G,\bar k(t)^*) @>v_P>>          H^2(G,\D{Z}_P)\\[-1mm]
         @V\cong VV                                             @VV\cong V\\[-1mm]
        \Br(k(t))                   @>\partial_A>>\Char(k(\alpha)).
\end{CD}
\]
Putting this together with sequence \eqref{eq.htwo.projective.line} completes the proof.
\end{proof}

\begin{prop}                                              \label{prop.bru.purity.invariance}
Let $K$ be a function field over the
field $k$. Let $t$ be an
indeterminate. The natural map from $\Br K$ to $\Br (K(t))$ induces an
isomorphism $\Bru (K/k)
\cong \Bru (K(t)/k)$. In particular, if $K/k$ is rational, \emph{i.e.}\ if $K$
is purely transcendental
over $k$, or more generally if $K/k$ is stably rational, {\rm i.e.} if
$K(t_1,\dots,t_r)$ is purely
transcendental over $k$
for some suitable independent variables
$t_1,\dots,t_r$, then the natural map
from $\Br (k)$ to $\Bru (K/k)$ is an isomorphism.
\end{prop}

\begin{proof}
It is enough to prove the first assertion. Let
$L=K(t)$. According to Lemma \ref{lem.bru.purity.invariance} 
the map $\Br K \to\Br (K(t))$ is an injection and according to
Lemma \ref{lem.functoriality.bru} it induces an
inclusion $\Bru (K/k) \subset\Bru(L/k)$. By the very definition of the
unramified Brauer group,
any $\alpha$ in $\Bru (L/k)$ clearly belongs to $\Bru (L/K)= \Bru
(K(t)/K)$, and Lemma \ref{lem.bru.purity.invariance} says that
this last group coincides with $\Br K$. Thus we only have to prove that
if $\alpha\in\Br K$
becomes unramified over
$k$ 
when viewed in $\Br (K(t))$, then it already
belongs to $\Bru (K/k)$.
Let $A$ be a discrete rank one valuation ring with field of fractions $K$.
The localization of
$A[t]$ at the prime ideal of $A[t]$ spanned by a uniformizing parameter
$\pi$ of $A$ is a discrete
valuation ring $B$ with uniformizing parameter $\pi$, and the induced map
on residue class fields
$\kappa_A \to\kappa_B$ may be identified with the inclusion $\kappa_A\to
\kappa_A(t)$. Hence in
diagram \eqref{eq.functoriality.residue} applied to the present situation, the R.H.S.\ vertical map,
i.e.,
$\Res_{\kappa_A/\kappa_B}\colon\Char(\kappa_A)\to  \Char(\kappa_B) = \Char(\kappa_A(t))$ is
an injection (this
applies more generally as soon as $e =1$ and $\kappa_A$ is algebraically
closed in $\kappa_B)$. Thus
for $a\in\Br K$, the equality $\partial_B(\alpha_L)=0$ implies
$\partial_A(\alpha)=0$. Since $A$
was an arbitrary rank one discrete valuation ring in $K$, the conclusion
follows.
\end{proof}

\begin{rem} 
More generally, one may show that if $K$ is the function field
of a $k$-variety which is retract rational (over $k$), then   
 the natural map
from $\Br (k)$ to $\Bru (K/k)$ is an isomorphism.
The proof requires functorial properties
of the unramified Brauer group  more elaborate than the one
given in Lemma \ref{lem.functoriality.bru}.  For this, we
refer to \cite{C}.
\end{rem}

Up till now, we have sticked to a down-to-earth definition of the
unramified Brauer group. In some
circumstances it is necessary to use a high-brow definition.

\begin{Def}[Grothendieck \cite{Gr2}]                             \label{def.groth.brauer.group}
The Brauer group $\Br  X$ of
a scheme $X$ is the
second \'etale cohomology group $H_{\text{\'et}}^2(X,\D{G}_m)$.
\end{Def}

\begin{rem}                                                                     \label{rem.groth.brauer.group}
Note that $H_{\text{\'et}}^2(X,\D{G}_m)$ is called ``cohomological Brauer group''
and denoted $\Br'X$
in \cite{Gr2}, whereas $\Br X$ denotes the ``Brauer group'' defined in \cite{Gr1} as
classes of Azumaya algebras over $X$.
\end{rem}

Using Grothendieck's theorems \cites{Gr1,Gr2}, and in particular the ``purity
theorem'' \citelist{\cite{Gr3}*{\S6}\cite{Del}*{p.~63, \S V, Th\'eor\`eme (3.4)}}, 
one may show (see \cite{C}*{\S 3.4 and Prop.~4.2.3}):

\begin{thm}                                                          \label{thm.groth.brauer.unramified}
Let $X$ be a smooth connected variety over
a field $k$ of
characteristic zero. Let $k(X)$ be the function field of $X$. 
\begin{enumerate}[\rm(i)]
\item
There is a
natural inclusion $\Br  X\subset\Br  k(X)$ --- in particular $\Br  X$ is a torsion group.
\item
The subgroup $\Bru(k(X)/k)\subset\Br k(X)$ lies in $\Br X$.
\item
 If
moreover $X$ is proper, $\Br  X=\Bru(k(X)/k)$.
\end{enumerate}
\end{thm}

\begin{rem*}
For examples of computations of $\Bru(k(X)/k)$ over a non-algebraically closed field $k$, 
see \cites{BoroK,CKun1,CKun2}.
\end{rem*}

For later use, let us recall that given any smooth variety $X$ over a field $k$
and any integer $n$ prime
to $\car k$, the Kummer sequence in \'etale cohomology induces a short exact
sequence:
\[
0 \to\Pic X/n \to  H_{\text{\'et}}^2(X,\mu_n)
\to {}_n\!\Br X \to  0,
\]
where $\Pic X$ denotes the Picard group of $X$, $\mu_n$ the group of $n$-th
roots of unity and
${}_n\!\Br X$ the $n$-torsion subgroup of $\Br X$.

\section{A general formula}          \label{sec.bru.formula}

The results of this section are mainly due to F.A.~Bogomolov \cites{Bo3, Bo4}.
Here $k$ denotes an algebraically closed field of characteristic zero.

We shall first consider the case of a finite group $G$.

\begin{thm}                                                                \label{thm.bru.bicyclic.finite}
Let $G$ be a finite group of automorphisms
of a function field $L$
over the algebraically closed field $k$ of characteristic zero. One then has:
\[
\Bru L^G= \accol{\alpha\in\Br L^G\mid\alpha_H\in\Bru L^H \text{ for all } H\in\C{B}_G},
\]
where 
$\C{B}_G$ denotes the set of finite bicyclic subgroups of $G $ (a bicyclic group
is a group spanned by at most two elements) and 
 $\alpha_H$ denotes the restriction of $\alpha\in\Br L^G$ to $\Br L^H$.
\end{thm}

\begin{proof}
Let $K=L^G$ and let $\alpha\in\Br K$ be such that
$\partial_A(\alpha)\neq 0$ for
some discrete rank one valuation ring $A\subset  K$ with fraction field
$K$. We must show that there
exists a subgroup $H\in\C{B}_G$ such that
\[
\alpha_H\notin \Bru L^H.
\]
The following facts may be read off from
Serre (\cite{Se1}*{I, \S 7}). Let $\F{p}$ be a prime ideal in
the semi-local Dedekind
ring  $\tilde A$ which is the integral closure of $A$ in $L$, and let
$D\subset  G$ be the associated
decomposition group, and let $I\subset D$ be the inertia group, which is a
normal subgroup of $G$.
The localization $B=\tilde A_{\F{p}} \subset  L$ is a discrete valuation
ring. There is a tower of
fields: $K\subset  L^D\subset  L^I\subset  L$ and a corresponding tower of
discrete valuation rings
obtained by taking the traces $A=B^G\subset B^D\subset B^I$ of $B$ on the
subfields. The
corresponding residue field extensions read: $F=F \subset   E=E$, and
we have $D/I= \Gal(E/F) = \Gal(L^I/L^D)$. The Galois extension $L^I/K$ is unramified, 
i.e.\ a uniformizing
parameter of $A$ is
still a uniformizing parameter in $B^I$.

Moreover, since the residue characteristic is zero, the inertia group $I$
may be
identified with a cyclic group, namely a group $\mu$ of roots of
unity in $F$ (\cite{Se1}*{IV, \S 2, Corollaires 1 and 2}). Furthermore, the conjugacy action of $D$ on the normal
subgroup $I$ is then
trivial, since this action may be identified with the action of $D/I=\Gal
(E/F)$ on $\mu\subset  F$,
and all roots of unity are in $k\subset  F$. Thus $I$ is central in $D$.

If $\alpha_I\notin \Bru L^I$, we are done, since $I$ is a cyclic subgroup of $G$. We may thus
assume that $\alpha_I\in
\Bru (L^I)$. Since $B^D/A$ is an unramified extension of discrete valuation
rings which induces an
isomorphism on the residue class fields, the assumption
$\partial_A(\alpha)\neq 0$ implies
$\partial_{B^D}(\alpha)\neq 0\in \Char(F)$. On the other hand,
$\partial_{B^I}(a)=0\in \Char(E)$. Since $B^I/B^D$
is unramified, the commutative diagram:
\[
\begin{CD}
        \Br K^I   @>\partial_{Bë}>>              \Char(E) \\
        @AAA                           @AA\Res_{F/E}A\\
       \Br K^D  @>\partial_{B^D}>> \Char(F)
\end{CD}
\]
implies that $\partial_{B^D}(\alpha)$ may be identified with a nontrivial
character of $D/I=\Gal(E/F)$. Let $g\in D$ be an element of $D$ whose class $\bar g$ in $D/I$
satisfies
$\partial_{B^D}(\alpha)(\bar g) \neq  0\in\D{Q}/\D{Z}$, let $H=\scal{I,g}\subset  D$ be the
subgroup spanned by $I$
and $g$, and let $F_1$ be the residue class field of $B^H$. Inserting $
\Br (K^H)\to \Char(F_1)$
in the above diagram, one immediately sees that $\partial(\alpha_H) \neq  0$,
since $\partial(\alpha_H)$
may be identified with a character of $\Gal(E/F_1)=D/H$ which does not
vanish on $\bar g$. This is
enough to conclude, since $H$ is an extension of the cyclic group $\scal{\bar
g}$ by the central cyclic
subgroup $I$ (see above), hence is an abelian group spanned by two
elements.
\end{proof}

We now wish to extend Theorem \ref{thm.bru.bicyclic.finite} to almost free actions of reductive
algebraic groups.

\begin{lem}                                      \label{lem.rational.section.bru}
Let $p\colon X\to Y$ be a 
dominant
morphism of smooth
integral varieties. If this
morphism admits a section over a non-empty open set of $Y$, then this section
induces a cartesian diagram, 
where the horizontal maps are injective:
\[
\xymatrix{
{\Br Y} \ar@{^{(}->}[r] &{\Br X}\\
{\Bru k(Y)}\ar@{^{(}->}[u] \ar@{^{(}->}[r]&{\Bru k(X).}\ar@{^{(}->}[u]
}
\]
\end{lem}

\begin{proof}
Let the section $s$ be defined over the open set
$V\subset  Y$, and let
$U=p^{-1}(V)\subset  X$. Because of the functorial behaviour of the Brauer
group and of the
unramified Brauer group we have a commutative diagram:
\[
\xymatrix{
{\Br V} \ar[r]^{p^*}                &{\Br U}\ar[r]^{s^*}               &{\Br V}\\
{\Br Y} \ar@{^{(}->}[u]
             \ar[r]^{p^*}               &{\Br X}\ar@{^{(}->}[u]         &\\
{\Bru k(Y)}\ar@{^{(}->}[u] 
                 \ar[r]^{p^*}           &{\Bru k(X)}\ar@{^{(}->}[u]
                                                                   \ar[r]^{s^*}         &{\Bru k(Y)}\ar@{^{(}->}[uu]
}
\]
where the top composite map is identity and where all vertical
maps are inclusions.  

The right hand part of the diagram requires an explanation,
the notation $s^*$ being a slight abuse of language.
Since the characteristic of $k$ is zero, Hironaka's theorem
guarantees that the morphism $p\colon X\to Y$
extends to a morphism $p^{\text{c}}\colon X^{\text{c}}\to Y^{\text{c}}$, where
$X^{\text{c}}$, resp.~$Y^{\text{c}}$, is a smooth, projective integral variety
containing $X$, resp.~$Y$, as a dense open set. Since
$p^{\text{c}}$ is a proper morphism and $Y^{\text{c}}$ is smooth, the section
$s\colon V \to X$ extends to a section $s^{\text{c}}\colon W \to X^{\text{c}}$
of $p^{\text{c}}$, where $W \subset Y^{\text{c}}$ is an open set which contains
all codimension $1$ points of $Y^{\text{c}}$:
\[
\xymatrix{
                                &U\ar@{^{(}->}[d]                      &V\ar[l]_{s}\ar@{=}[d]      & \\
Y\ar@{^{(}->}[d]   &X\ar[l]_{p}\ar@{^{(}->}[d]        &V\ar[l]_{s}\ar@{^{(}->}[r] \ar@{^{(}->}[d]      &Y \ar@{^{(}->}[d] \\
Y^{\text{c}}           &X^{\text{c}}\ar[l]_{p^{\text{c}}} &W\ar[l]_{s^{\text{c}}}\ar@{^{(}->}[r]               &Y^{\text{c}}.
}
\]
By purity of the Brauer group (Theorem \ref{thm.groth.brauer.unramified}), the restriction map
$\Br Y^{\text{c}}\to\Br W$ is an  isomorphism, and both groups
coincide with $\Bru k(Y)$. We thus have the commutative diagram:
\[
\xymatrix{
&{\Br U}\ar[r]^{s^*}                                 &{\Br V}                                        &&\\
&{\Br X}\ar[r]^{s^*}\ar@{^{(}->}[u]        &{\Br V}\ar@{=}[u]                      &&\\
{\Bru k(X)}\ar@{=}[r] 
&{\Br X^{\text{c}}}\ar[r]^{s^{\text{c}*}}\ar@{^{(}->}[u]&{\Br W}\ar@{^{(}->}[u] &{\Br Y^{\text{c}}}\ar[l]_{\cong}
&{\Bru k(Y)}.\ar@{=}[l]
}
\]
The composite map $\Br X^{\text{c}}\to\Br X \to\Br U\to\Br V$
coincides with the composite map $\Br X^{\text{c}}\to\Br W\to\Br V $,
which may be rewritten $\Bru k(X) \to\Bru k(Y) \to\Br V$.

Now if $\alpha\in\Br Y$ is such that $p^*(\alpha)\in\Br X$ actually lies
in $\Bru k(X)$, the above
diagram shows that $\Res_{Y/V}(\alpha) = s^*p^*(\Res_{Y/V}(\alpha))$
belongs to $\Bru k(Y)$, hence
also $\alpha\in  \Bru k(Y)$.
\end{proof}

\begin{lem}[Bogomolov \cite{Bo3}]                  \label{lem.bogomolov.finite.reduction}
Let $k$ be an algebraically closed field, $\car k=0$,
and let $G$ be a reductive algebraic
group over $k$ which is an
extension of a finite group $W$ by a torus $T$. Let $X$ be an integral
affine $k$-variety with an action of $G$. Assume that all stabilizers are trivial. 
Then, there exists a finite subgroup $A$ of $G$ such that the
natural map $X/A \to  X/G$ has
a section over a non-empty open set of $X/G$.
\end{lem}

\begin{proof}
First af all, for any algebraic subgroup $H$ of $G$, the natural map $X\to X\quota H$
makes $X$ into a $H$-torsor over $X\quota H$. Hence $X/H$ exists, and $X/H=X\quota H$.
By assumption, the group $G$ defines an extension:
\[
1\to  T\to  G\to  W\to  1,
\]
hence an algebraic action of $W$ on $T$, and a class in $H^2(W,T)$. 
Let ${}_nT$ be the subgroup of $n$-torsion
points of $T$ for
$n\geq 1$, and let $T_{\text{tors}}= \bigcup_{n\geq 1}{}_nT$ be the whole
torsion subgroup of $T$. The
exact sequence
\[
1 \to  T_{\text{tors}} \to  T \to  T\otimes_{\D{Z}}\D{Q} \to  1
\]
and the vanishing of $H^i(W,T\otimes_{\D{Z}}\D{Q})$ for $i\geq 1$ shows that this
class comes from a (unique)
class in $H^2(W,T_{\text{tors}})$, hence from some class in $H^2(W,{}_nT)$ for
some $n$. Thus there is a
finite group $H_n$ and a commutative diagram of extensions:
\[
\begin{CD}
         1         @>>>     T           @>>>    G           @>>>   W      @>>>   1\\
         @.                    @AAA                 @AAA                @|                 @.  \\
         1         @>>>     {}_nT    @>>>    H_n      @>>>   W      @>>>   1.
\end{CD}
\]
This diagram gives rise to the following fibre product:
\[
\begin{CD}
         X/{}_nT          @>>>    X/H_n\\
         @Vq_nVV                   @VVp_nV\\
         X/T                 @>>>     X/G.
\end{CD}
\]
Moreover, the map $q_n$ makes $X/{}_nT$ into a torsor over $X/T$ under the
torus $T=T/{}_nT$. Note
that the group $W$ acts upon $X/{}_nT$ and $X/T$ (indeed, both horizontal
maps in the above diagram
are Galois coverings with group $W$). The map $q_n$ is $W$-equivariant.
Since torsors under tori are
locally trivial for the Zariski topology, the set of rational sections of
$q_n$ is not empty. The
group $W$ acts upon this set. 

If we may find such a section which is
$W$-invariant (for the obvious
action induced by the action of $W$ on both spaces $X/{}_nT$ and $X/T$),
then this section descends
to a section of $p_n$ and we get the conclusion of the lemma. 

We shall show
that at the cost of
changing $n$ into $nm$ for a suitable $m\geq 1$, there exists such a
section. Let $E=k(X)^T=k(X/T)$
and let $F=k(X)^G=k(X/G)$. The extension $E/F$ is a Galois extension with
Galois group $W$. Thus to
the action of $W$ on the $k$-torus $T$ we may associate a unique (twisted)
$F$-torus $R$ which
becomes isomorphic to the torus $T_E=T\times_kE$ over $E$. Namely, we take
\[
R=\Spec F[\Char(T)]^W,
\]
where $F[\Char(T)]$ is the group algebra
over the character group $\Char(T)$ of $T$ (which is a free abelian group of finite
type), and the $W$-action on
$F[\Char(T)]$ is simultaneous on $F$ and $\Char(T)$. The generic fibre of $q_n$ is
a principal homogenous
space under the action of $T_E$. Since $W$ acts equivariantly on the whole
situation, the generic fibre
$P_n$ of $p_n$, which is an $F$-variety, inherits a structure of principal
homogeneous space under
the $F$-torus $R$. The isomorphy class of this principal homogeneous space
in the \'etale (= Galois)
cohomology group $H^1(F,R)$ is killed by some integer $m\geq 1$. Let us
then consider the
commutative diagram of exact sequences:
\[
\begin{CD}
         1         @>>>   {}_{mn}R  @>>>    R       @>mn>>   R              @>>>   1\\
         @.                    @AAA                     @|                         @AAmA                @.  \\
         1         @>>>     {}_nR     @>>>     R      @>n>>       R              @>>>    1,
\end{CD}
\]
which reflects the $W$-action on the commutative diagram
\[
\begin{CD}
         1         @>>>   {}_{mn}T  @>>>    T       @>mn>>   T              @>>>   1\\
         @.                    @AAA                     @|                         @AAmA                @.  \\
         1         @>>>     {}_nT     @>>>     T      @>n>>       T              @>>>    1.
\end{CD}
\]
We also have the fibre product:
\[
\begin{CD}
         X/{}_{mn}T           @>>>    X/H_{mn} \\
         @Vq_{mn}VV                    @VVp_{mn}V\\
         X/T                        @>>>     X/G,
\end{CD}
\]
where $H_{nm}$ is the finite group obtained from $H_n$ by pushing out
through ${}_nT\to {}_{mn}T$,
thus yielding the commutative diagram:
\[
\begin{CD}
         1         @>>>   T                @>>>    G             @>>>      W       @>>>    1\\
         @.                    @AAA                     @AAA                     @|                    @.  \\
         1         @>>>  {}_{mn}T   @>>>     H_{mn} @>>>       W       @>>>    1\\
         @.                    @AAA                     @|                            @|                    @.  \\
         1         @>>>   {}_nT        @>>>     R            @>>>       R       @>>>    1.
\end{CD}
\]
Both the generic fibre $P_n$ of $q_n$ and $P_{nm}$ of $q_{nm}$ are
principal homogeneous
spaces under the same $F$-torus $R$. Indeed, on $X/{}_nT$, the torus which
acts is the quotient
$T=T/{}_nT$, and on $X/{}_nT$, it is $T=T/{}_{nm}T$. Now as the previous
diagrams reveal, $r\colon
X/{}_nT \to  X/{}_{nm}T$ is the obvious projection, $t$ is in $T$ and $x$
is in $X/{}_nT$, then
$r(t.x)=t^m\cdot r(x)$. Thus the natural map from $X/H_n \to  X/H_{nm}$
induces on generic fibres a
map $r_1: P_n\to P_{nm}$ which satisfies
\[
r_1(t\cdot x)=t^m\cdot r_1(x) \text{ for } t\in R.
\]
This implies that the class of $P_{nm}$ in $H^1(F,R)$ is $m$
times the class of $P_n$. Since the class of $P_n$ is killed by $m$, we
conclude that the class of
$P_{nm}$ is trivial, which completes the proof.
\end{proof}

\begin{thm}[Bogomolov \cite{Bo4}*{Theorem 2.1}]                \label{thm.bru.bicyclic.reductive}
Let $k$ be an algebraically closed field of characteristic zero, 
let $G$ be a reductive group over $k$, and let
$X$ be an integral affine
$k$-variety with a $G$-action. Assume that all stabilizers are trivial. One then has
\[
\Bru k(X)^G = \accol{\alpha\in\Br k(X)^G\mid\alpha_A\in\Bru k(X)^A
\text{ for all } A\in\C{B}_G},
\] 
where $\C{B}_G$ denotes the set of finite bicyclic subgroups of $G(k)$
and $\alpha_A$ denotes the restriction of 
$\alpha\in\Br k(X)^G$ to $\Br k(X)^A$.
\end{thm}

\begin{proof}
Let $G^{\circ}$ be the connected component of identity,
and let $\F{g}$ be its
Lie algebra. Let the group $G$ act on $\F{g}$ via the adjoint
representation and let $G$ act upon
the product $X\times\F{g}$ via the diagonal action. Let $\F{t}\subset
\F{g}$ denote a fixed Cartan
subalgebra, and let $N\subset  G$ be the normalizer of $\F{t}$ for the
adjoint action. The group $N$
is a reductive group whose identity component is a (maximal) torus
$T\subset  G$. It is a classical
fact that any regular semisimple element in $\F{g}$ is conjugate under
$G$ to an element of
$\F{t}$. Moreover, on the open set $\F{t}^{\circ}\subset  \F{t}$
consisting of regular semisimple
elements, the following property holds: if $g\in G$, $x, y\in\F{t}^{\circ}$
and $g\cdot x=y$,
then $g\in N$. Thus the closed subvariety $X\times\F{t}\subset
X\times\F{g}$ is a $(G,N)$-slice,
and by the slice lemma \ref{thm.slice.lemma} we have $k(X\times\F{t})^N\cong k(X\times\F{g})^G$. On the other
hand, the no-name lemma \ref{cor.noname.lemma}
implies that $k(X\times\F{t})^N$, resp.\ $k(X\times\F{g})^G$, is purely
transcendental over
$k(X)^N$, resp.\ $k(X)^G$. According to Proposition \ref{prop.bru.purity.invariance}, we thus have: $\Bru
k(X)^N\cong\Bru k(X\times\F{t})^N$
and $\Bru k(X)^G\cong \Bru k(X\times \F{g})^G$. The functoriality of the
unramified Brauer group
now implies that the natural inclusion $k(X)^G\subset  k(X)^N$ induces an
isomorphism $\Bru k(X)^N\cong \Bru k(X)^G$. According to Lemma \ref{lem.bogomolov.finite.reduction}, 
there exists a finite
group $H\subset  N$ such
that the map between quotient spaces: $X/H\to X/N$ has a rational section.
Now Lemma \ref{lem.rational.section.bru} shows that
an element $\alpha\in\Br (k(X)^N)$ is unramified if and only if its
restriction $\alpha_H$ is
unramified. By Theorem \ref{thm.bru.bicyclic.finite}, $\alpha_H$ itself is unramified if and only if
$\alpha_A$ is unramified
for all bicyclic subgroups $A$ of $H$. The theorem follows.
\end{proof}

\section{Linear action of a finite group}       
                                                         \label{sec.linear.action}

In this section, $k$ denotes an algebraically closed field of characteristic zero.

\begin{thm}[Bogomolov \cite{Bo3}]   \label{thm.bru.formula.finite}
Let $G\subset
\GL(V)$ be a finite group of
automorphisms of a finite dimensional $k$-vector space $V$. One then has:
\begin{align*}
\Bru k(V)^G 
&\cong	
\ker\Bigl\lbrack H^2(G,k^*)\xrightarrow{\Res}\prod_{A\in\C{B}_G}H^2(A,k^*)\Bigr\rbrack\\
&\cong	
\ker\Bigl\lbrack H^2(G,\D{Q}/\D{Z})\xrightarrow{\Res}\prod_{A\in\C{B}_G} H^2(A,\D{Q}/\D{Z})\Bigr\rbrack\\
&\cong
\ker\Bigl\lbrack H^3(G,\D{Z})\xrightarrow{\Res}\prod_{A\in\C{B}_G} H^3(A,\D{Z})\Bigr\rbrack,
\end{align*}
where $\C{B}_G$ denotes the set of all bicyclic subgroups of $G$ (in the
last two formulas, $G$ acts
trivially upon $\D{Z}$ and $\D{Q}/\D{Z}$).

The same formulas hold if one replaces the
set $\C{B}_G$ of bicyclic subgroups by the set $\C{A}_G$ of all abelian
subgroups of $G$.
\end{thm}

\begin{proof}
If $A\subset  G$ is any abelian group, the field
$k(V)^A$ is pure over $k$ by
Fischer's theorem (see Proposition \ref{prop.torus.purity}), hence $\Bru k(V)^A = 0$. On the other hand, given any
subgroup $H\subset  G$, the
restriction-inflation sequence for the faithful action of $H$ on $k(V)^*$
yields an exact sequence:
\[
0 \to  H^2(H,k(V)^*) \to\Br  k(V)^H \to\Br  k(V).
\]
Using the functoriality of the unramified Brauer group 
\ref{lem.functoriality.bru}  and 
Theorem  \ref{thm.bru.bicyclic.finite},
one gets the formula:
\begin{equation}                                          \label{eq.bru.formula.finite}
\Bru k(V)^G \cong	
\ker\Bigl\lbrack H^2(G,k(V)^*)\xrightarrow{\Res} \prod_{A\in\C{B}_G}
H^2(A,k(V)^*)\Bigr\rbrack.
\end{equation}
Since $k\croch{V}$ is a UFD and $k\croch{V}^*=k^*$, 
one has a short exact sequence of $G$-modules:
\[
1 \to  k^* \to  k(V)^* \to\Div V \to  0,
\]
where the $G$-module $\Div V$ is a direct sum of permutation modules
$\D{Z}\croch{G/H}$ for
various subgroups $H$. Such a module satisfies the two properties:
\begin{align*}
&H^1(G,M)=0,\\
&\ker\Bigl\lbrack H^2(G,M) \to\prod_{g\in G}H^2(\scal{g},M)\Bigr\rbrack = 0,
\end{align*}
as one easily checks by reducing to the case $M=\D{Z}$ with trivial action.
The first formula now easily
follows. The second formula is obtained by identifying the group $\mu$ of
roots of unity in $k^*$
with $\D{Q}/\D{Z}$, and using the unique divisibility, hence cohomological
triviality of $k^*/\mu$. As for
the third formula, it is obtained by shifting via the exact sequence $0\to
\D{Z}\to\D{Q}\to\D{Q}/\D{Z}\to 0$.

The last statement 
now follows from the vanishing $\Bru k(V)^A =0$ for any abelian subgroup $A$ of $G$.
\end{proof}

\subsection{The case of a nilpotent 
             group of class 2}          \label{ssec.nilpotent.class.two}

Bogomolov \cite{Bo3} made a thorough application of the previous theorem
when $G$ is a nilpotent group,
particularly of class $2$. Such a group is a central extension
\[
1 \to  C \to  G \to\Gamma \to  1
\]
of a finite abelian group $\Gamma$ by another finite abelian group 
$C$ and classes of such extensions
are classified by the group
$H^2(\Gamma,C)$ where $C$ is viewed as a trivial $\Gamma$-module. We denote by
$\croch{G}\in H^2(\Gamma,C)$ the class
of $G$.

It is well known \cite{Bro}*{\S V.6, Theorem 6.4}
that for an abelian group $\Gamma$ we have
$H_2(\Gamma,\D{Z})=\Lambda^2\Gamma$. Then, for any $\Gamma$-module $M$, the
universal coefficient sequence \cite{Bro}*{\S V.6, exercise 5}
yields
\[
0 \to\Ext^1(\Gamma,M) \xrightarrow{} H^2(\Gamma,M)\xrightarrow{\omega_M}
\Hom(\Lambda^2\Gamma,M) \to  0.
\]
Then for $M = \D{Q}/\D{Z}$ with trivial action,  the map
\[
\omega_{\D{Q}/\D{Z}}\colon
H^2(\Gamma,\D{Q}/\D{Z})\to\Hom(\Lambda^2\Gamma,\D{Q}/\D{Z})
\]
is an isomorphism. For $M = C$ with trivial action, the map 
\[
\omega_C\colon
H^2(\Gamma,C) \xrightarrow{}
\Hom(\Lambda^2\Gamma,C)
\] sends the class of $\croch{G}$ to
\[
\lambda_G\colon\Lambda^2\Gamma \to \croch{G,G}\subset C
\]
defined by
\[
\lambda_G(\gamma_1\wedge \gamma_2) = \croch{g_1,g_2}
\]
where $g_1, g_2\in G$ are lifts of $\gamma_1, \gamma_2\in\Gamma$.
From the definition $\im\lambda_G=\croch{G,G}$. So $\lambda_G$ is surjective if
and only if
$\croch{G,G}=C$.

\begin{Def}                                                    \label{def.sbic.central.extension}
Let $G$ be a  central extension of a finite abelian
group $\Gamma$ by another one $C$. To the canonical
homomorphism $\lambda_G\colon\Lambda^2\Gamma\to C$ one may attach the following subgroups
\[
S_{\text{bic}} = S_{G,\text{bic}}\subset S_G:= \ker\lambda_G \subset
\Lambda^2\Gamma
\]
where $S_{\text{bic}}$ is the subgroup of $S_G$ generated by all the
$\gamma_1\wedge\gamma_2$ which
belong to $S_G$.
\end{Def}

In other words, $S_G$ is defined by the exact sequence
 \[
0\to S_G\to\Lambda^2\Gamma\xrightarrow{\lambda_G}C
\]
and $S_{\text{bic}}$ is the subgroup generated by the $S_A$'s for all bicyclic
subgroups $A$ of $G$.

\begin{thm}[Bogomolov \cite{Bo3}]                                    \label{thm.bru.formula.central.extension}
Let $G$ be a  central extension of a finite abelian
group $\Gamma$ by another one $C$. Then
\[
\Bru k(V)^G = \ker(\widehat{S_G}\to\widehat{S_{\textup{bic}}}) =
\widehat{(S_G/S_{\textup{bic}})}
\]
where $S_G$ is the kernel of the canonical morphism $\lambda_G\colon\Lambda^2\Gamma\to C$,
and
\[
S_{\textup{bic}}\subset S_G \subset \Lambda^2\Gamma
\]
is the subgroup generated by the $S_A$'s for all bicyclic
subgroups $A$ of $G$.
\end{thm}

\begin{proof}
The lower terms exact sequence of the Hochschild-Serre  spectral sequence
\[
H^p(\Gamma,H^q(C,\D{Q}/\D{Z})) \implies H^{p+q}(G,\D{Q}/\D{Z})
\]
yields the top exact sequence in the following diagram:
\begin{small}%
\[
\begin{CD}
\widehat C  @>>> H^2(\Gamma,\D{Q}/\D{Z}) @>>>  
                                                           \ker[H^2(G,\D{Q}/\D{Z})\to H^2(C,\D{Q}/\D{Z})]  @>>>  H^1(\Gamma,\widehat C)\\
         @|             @V\omega_{\D{Q}/\D{Z}}V\cong V                     
                                                                                                                      @AAA                           @.  \\
\widehat C  @>>> \Hom(\Lambda^2\Gamma,\D{Q}/\D{Z}) @>>> \Hom(S_G,\D{Q}/\D{Z})=\widehat {S_G}  @>>> 0.
\end{CD}
\]
\end{small}%
This gives the following exact sequence:
\[
0\to\widehat {S_G}\to\ker\bigl\lbrack H^2(G,\D{Q}/\D{Z})\to H^2(C,\D{Q}/\D{Z})\bigr\rbrack\to
H^1(\Gamma,\widehat C).
\]
We want to calculate
\[
B_G:= \Bru \D{C}(V)^G \subset H^2(G,\D{Q}/\D{Z}).
\]
This inclusion is functorial and covariant in $G$. Fischer's theorem yields
$B_A =
0$ for any abelian group $A$. In particular, $B_C = 0$. So by restriction
$B_G$ maps to $0$
in $H^2(C,\D{Q}/\D{Z})$:
\[
B_G \subset \ker\bigl\lbrack H^2(G,\D{Q}/\D{Z}) \to  H^2(C,\D{Q}/\D{Z})\bigr\rbrack.
\]
Then we prove that $B_G$ maps to $0$ in $H^1(\Gamma,\widehat
C)=\Hom(\Gamma,\widehat C)$. By
functoriality, it is sufficient to prove this result by restriction to any
cyclic subgroup
$\Gamma'$ of $\Gamma$. Let $G'$ the restriction of the extension $G$ to
$\Gamma'$.
A central extension of a cyclic group is an abelian group. Thus $\Gamma'$
is abelian and again
by Fisher's theorem $B_{G'}=0$.

We finally get the inclusion
\[
B_G\subset \widehat {S_G}.
\]
Theorem \ref{thm.bru.formula.finite} says that
\[
B_G\cong	\ker\Bigl\lbrack H^2(G,\D{Q}/\D{Z}) \xrightarrow{} \prod_{G'\in\C{B}_G} H^2(G',\D{Q}/\D{Z})\Bigr\rbrack
\]
where $G'$ is any bicyclic subgroup of $G$. Then by functoriality $B_G$ is
the subgroup of $\widehat
{S_G}$ which maps to $0$ in each of the $H^2(G',\D{Q}/\D{Z})$'s, hence in each of
the $S_{G'}$'s for any
bicyclic subgroup $G'$ of $G$. Let $\Gamma'$ be the image in $\Gamma$ of
such a subgroup $G'$. We
have a commutative diagram
\[
\begin{CD}
         0         @>>>   S_G          @>>>   \Lambda^2\Gamma   @>\lambda_G>>        C\\
         @.                    @AAA                     @AAA                                                                @|  \\
         0         @>>>    S_{G'}     @>>>  \Lambda^2\Gamma'   @>\lambda_{G'}\,=\,0>> C.
\end{CD}
\]
Let $g_1, g_2\in G$ and $\gamma_1, \gamma_2\in\Gamma$ their images. Then
the two conditions
\begin{enumerate}[(i)]
\item
$g_1g_2=g_2g_1$,
\item
$\gamma_1\wedge\gamma_2\in S_G$,
\end{enumerate}
are equivalent.

Hence the $S_{G'}$'s for any bicyclic subgroup $G'\subset G$ are exactly
the cyclic
subgroups of $S_G$ generated by elements of the form $\gamma_1\wedge
\gamma_2$.
\end{proof}

\begin{cor}                                                    \label{cor.sbic.central.extension}
Let $\Gamma$ be a finite abelian
group. There is a bijection
between the classes of central extensions $G$ of $\Gamma$ such that $B_G\neq 0$ and the subgroups
\[
S\subset\Lambda^2\Gamma
\]
such that
\[
S_{\textup{bic}}\neq S.
\]
\end{cor}

\begin{exs}                                                   \label{exs.sbic.central.extension}
Let us find simple examples of 
finite groups $G$ such that $B_G\neq 0$. We
consider only nilpotent
groups $G$ of class $2$.

We are searching $G$ among central extensions of a fixed finite abelian
group $\Gamma$ by
using Corollary \ref{cor.sbic.central.extension}.

Consider the case $\Gamma = (\D{Z}/p)^4$. We find
$\Lambda^2\Gamma = \D{F}_p^6$. It is simpler to work in the projective space
$\D{P}(\Lambda^2\Gamma) =
\D{P}^5(\D{F}_p)$. For $S\subset\Lambda^2\Gamma = \D{F}_p^6$ the following
conditions are equivalent:
\begin{enumerate}[(i)]
\item
$S\neq S_{\text{bic}}$.
\item
$\D{P}(S)\cap Q$ does not generate $\D{P}(S)$ as a linear space, where $Q$ denotes the quadric
whose equation is $x_1x_2+x_3x_4+x_5x_6 =
0$ which consists of all undecomposable tensors $\gamma_1\wedge \gamma_2$.
\end{enumerate}
Here is the list of all possibilities for condition (ii):
\begin{enumerate}[1.]
\item
$\D{P}(S) =$ a point $\notin Q$ (this is Saltman's original example \cite{Sa3}; here
$\#G = p^9$).
\item
$\D{P}(S) =$ a line tangent to $Q$ (here $\#G = p^8$).
\item
$\D{P}(S) =$ a line which does not intersect $Q$ (on $\D{F}_p$) (here $\#G = p^8$).
\item
$\D{P}(S) =$ a $2$-plane which intersects $Q$ along a line (here $\#G = p^7$).
\item
$\D{P}(S) =$ a $2$-plane which intersects $Q$ in a point (here $\#G = p^7$).
\end{enumerate}
More generally, for  $\Gamma = (\D{Z}/p)^{2m}$, the subvariety of
$\D{P}(\Lambda^2\Gamma)$  which
consists of tensors in $\Lambda^2\Gamma$ of rank $< 2m$ is a hypersurface
of
degree $m$ which is the zero locus of the Pfaffian. Among different
examples, the simplest one is obtained by taking 
$\D{P}(S)$ to be a point outside this hypersurface. One may also consider the
subspace $S$ which consists
of all matrices
\[
\begin{pmatrix}M_{1,1}&M_{1,2}\\
                           M_{2,1}&M_{2,2}\end{pmatrix}
\]
whose $4$ blocks are $m\times m$ matrices of the following type:
\begin{align*}
&M_{1,1}=M_{2,2}=0\\
&M_{1,2} \text{ trigonal with diagonal entries }
\accol{\lambda,\dots,\lambda}.
\end{align*}
in which case $S_{\text{bic}}$ is the subspace $\lambda=0$.
\end{exs}

\begin{rem*}
In \cite{Bo3} Bogomolov also gives an example of a group $G$ of order $p^6$ with $B_G\neq 0$
and proves that $B_G=0$ for any $G$ of order $p^n$ for $n\leq 5$. As a matter of fact, for $n\leq 4$,
the fields $k(V)^G$ are stably rational (Chu and Kang \cite{ChKan}). The stable rationality for $n=5$ remains open.
Recent computations of $B_G$ for finite Chevalley groups of type $A_n$ may be found in \cite{BoMP}.
\end{rem*}

\section{Multiplicative action of a finite group}
                                          \label{sec.multiplicative.action}

Most results in the present section are due to D.~Saltman.
Here $k$ denotes an algebraically closed field of characteristic zero.

Let $G$ be a group. A $G$-lattice $M$ is a $\D{Z}$-free module $M$ of finite rank
equipped with a $G$-action. We denote $k\croch{M}$ the group algebra of $M$ over $k$.
If $M$ is of rank $r$, then $k\croch{M}\cong k\croch{t_1,t_1^{-1},\dots,t_r,t_r^{-1}}$.
The associated field of fractions is denoted
$k(M)\cong k(t_1,\dots,t_r)$.

\begin{Def}                                                                \label{def.multiplicative.action}
Given a $G$-lattice $M$ consider 
the action of $G$ on $k\croch{M}$ which is trivial on $k$ and coincides with the given one on $M$.
The induced action of $G$ on $k(M)$ is called the \emph{multiplicative action} of
$G$ associated to the $G$-lattice $M$.
\end{Def}

Such actions are sometimes referred to as ``purely monomial actions''.

\begin{Def}                    \label{def.twisted.multiplicative.action}
For a given $G$-lattice $M$ we may also consider the \emph{twisted
multiplicative action}
of $G$ associated to a given crossed homomorphism of the finite group $G$ in $\widehat
M=\Hom(M,k^*)$, i.e.\ to a
$1$-cocycle $\alpha\in Z^1(G,\widehat M)$. If we denote by $t^m$ the
canonical image of $m\in M$ in
$k\croch{M}$ the $\alpha$-twisted action of $G$ on the $k$-algebra $k\croch{M}$ is
given by:
\[
g\cdot_{\alpha}t^m = \alpha_{g^{-1}}(m)t^{g\cdot m}.
\]
with trivial action on the coefficients.
We denote by $k\croch{M}_\alpha$ this   twisted   $G$-$k$-algebra and by
$k(M)_\alpha$ the field
of fractions. Let $\alpha$ and $\alpha'$ be in the same cohomology class.
Then the two
$G$-$k$-algebras $k\croch{M}_\alpha$ and $k\croch{M}_{\alpha'}$ are isomorphic: if
$\alpha'/\alpha=d\beta$
where $\beta: M\to k^*$ then $t^m\mapsto \beta(m)t^m$ defines a morphism of
$k$-algebras
$``\beta"\colon k\croch{M}_\alpha\to k\croch{M}_{\alpha'}$ which is $G$-equivariant and
which is clearly an
isomorphism. Note that  the formula
\[
g\cdot_{\alpha}t^m = \alpha_{g}(g\cdot m)t^{g\cdot m}
\]
defines the  opposite twisted  action on $k\croch{M}$, since we have
\[
\alpha_{g}(g\cdot
m)=({}^{g^{-1}}\!\alpha_g)(m) \text{ and } 1=\alpha_e
=\alpha_{g^{-1}}({}^{g^{-1}}\!\alpha_g).
\]
\end{Def}

A good reason for studying such actions is provided by the following theorem, which extends
early work of Procesi \cite{Pro1} and Formanek \cite{Fo1} for $G=\PGL_n$:

\begin{thm}[Saltman \cite{Sa8}]                                  \label{thm.multiplicative.weyl.action}
Let $G$ be a reductive connected linear algebraic group over $k$.
Let $V$ be a finite-dimensional $k$-vector space with an almost free linear $G$-action.
Let $T\subset G$ be a maximal torus, $W=N_G(T)/T$ the Weyl group and $\Char(T)$ the character group of $T$, 
viewed as a $W$-lattice. Choose a surjective map $f\colon P\to\Char(T)$ with $P$ a permutation $W$-lattice.
Let $M$ be the kernel of $f$. 
Then $k(V)^G$ is stably isomorphic over $k$ to $k_{\alpha}(M)^W$ for suitable 
$\alpha\in Z^1(W,\Hom(M,k^*))$.
\end{thm}

For further work in this direction, see \cites{Be3,BeSa,Sa16,SaTi}.

\begin{lem}                                  \label{lem.multiplicative.action.extension}
Every twisted $G$-$k$-algebra 
$k\croch{M}_\alpha$
defines an extension of
$G$-modules
\[
1\to k^*\to k\croch{M}_\alpha^* \to M\to 0 
\]
whose class is precisely the class of $\alpha\in Z^1(G,\Hom(M,k^*))$.
\end{lem}

\begin{lem}                   \label{lem.multiplicative.action.permutation}
Any  
twisted
multiplicative action of a
group $G$ on a
permutation $G$-lattice $P$ extends to a $G$-linear action on
a vector space $V$ such that $k(M)^G=k(V)^G$.
\end{lem}

\begin{proof}
Let $\accol{e_x}$ a basis of the $\D{Z}$-lattice $P$
which is $G$-stable. Then $G$ acts linearly on $V=\oplus_xke_x$. As $k\croch{P}$
is obtained from
$k\croch{V}$ by inverting the $x$'s we have the natural $G$-equivariant
inclusions
\[
k\croch{V}\subset
k\croch{P}\subset k(V)=k(P)
\]
and the same fields of invariants
$k(P)^G=k(V)^G$.
\end{proof}

\begin{lem}                   \label{lem.multiplicative.action.linear}
Given a multiplicative  action of a
group $G$  on a lattice $M$ one may find a $G$-lattice $N$, a finite
group $\Gamma$ which is a semi-direct extension of $G$ by an abelian group
and a faithful linear action of $\Gamma$ on a finite dimensional vector space $V$
such that $k(M\oplus N)^G=k(V)^{\Gamma}$.
\end{lem}

\begin{proof}
Up to changing  $M$ by $M\oplus N$, where $N$ is an auxiliary
$G$-lattice, one may find an exact sequence of $G$-modules
\[
0\to M\xrightarrow{} P\xrightarrow{\chi} F\to 0,
\]
where $P$ is a permutation $G$-lattice and $F$ is finite. The extension
$k(P)/k(M)$ is Galois
with group $\widehat F = \Hom(F,\mu)$ (where $\mu \subset  k^*$
denotes the group of roots of unity in $k$).
The $G$-action on $k(M)$ extends
to $k(P)$. Then
the extension  $k(P)/k(M)^G$ is also Galois with Galois group $\Gamma$
which is the  semi-direct
product of $\widehat F$ with $G$:
\[
\Gamma = \widehat F\rtimes G.
\]
The group $\Gamma$ acts linearly on $k(P)$. This is clear
for $g\in G$. One can verify that each $\gamma\in \widehat F$
acts on  $P \subset k(P)$ through multiplication by roots of unity:
\[
\gamma\cdot p = \scal{\gamma,\chi (p)}p
\text{ for } \gamma \in\widehat F \text{ and } p \in  P.
\]
Thus $\Gamma$ acts linearly
on $k\otimes_{\D{Z}}P$.
\end{proof}

For $M=\D{Z}\croch{G}$ the following theorem reduces to Theorem \ref{thm.bru.formula.finite}.

\begin{thm}[Saltman \cite{Sa10}*{Theorem 12}]                  \label{thm.multiplicative.bru.formula}
Let $G$ be a finite
group, let $M$ be a faithful
$G$-lattice, and let $k(M)$ denote the field of fractions of the group
algebra $k\croch{M}$. One then
has:
\[
\Bru k(M)^G \cong	
\ker\Bigl\lbrack H^2(G,k^*\oplus M)\xrightarrow{\Res}\prod_{A\in\C{B}_G} H^2(A,k^*\oplus
M)\Bigr\rbrack.
\]
\end{thm}

\begin{proof}
The variety $X=
\Spec k\croch{M}$ is isomorphic to a product of copies of the multiplicative group
$\D{G}_m$. Hence it is
factorial ($\Pic X=0$). Note that $\Div X$ is always a direct sum of permutation modules.
Moreover, the group of
units $k\croch{M}^*$ may be identified with the group $k^*\oplus M$ (as
$G$-modules). Thus one has an
exact sequence of $G$-modules:
\[
1\to  k^*\oplus M \to  k(M)^* \to\Div X \to  0 .
\]
If we argue as in Theorem \ref{thm.bru.formula.finite}, the theorem will follow from the following propositions:
\renewcommand{\qed}{}
\end{proof}

\begin{prop}                 \label{prop.multiplicative.bru.cyclic}
Let $A$ be a
cyclic group, and let $M$ be an $A$-lattice. Then the field of invariants
$k(M)^A$ is retract rational
and $\Bru k(M)^A=0$.
\end{prop}

\begin{proof}
Replacing $A$ by its image in $\Aut M$ allows us to assume that $M$ is a
faithful $A$-module.
Let $0\to M\to P\to F\to 0$ be a flasque resolution of $M$. Here $P$ is
a permutation
$G$-lattice and $F$ is a flasque $A$-lattice, i.e.\ $H^1(H,\Hom(F,\D{Z}))=0$
for all subgroups $H$ of
$A$. Such resolutions always exist \cite{CS1}. If $A$ is cyclic, a basic result of
Endo and Miyata says that any
flasque $A$-lattice is a direct factor of a permutation $A$-lattice, see \cite{CS1}. Quite
generally, given such an
exact sequence as above, the field $k(P)^A/k(M)^A$ is the function field of
a principal homogeneous
space $E$ over $k(M)^A$ under the $k(M)^A$-torus whose character
group over $k(M)$ is the $A$-lattice $F$. Since $F$ here is a
direct factor of a permutation module, it follows from Hilbert's theorem 90
that any such principal
homogenous space is trivial. In particular, $E$ has a $k(M)^A$-rational
point, and by Lemma \ref{lem.rational.section.bru} $\Bru
k(M)^A$ injects into $\Bru k(P)^A$. This last group is trivial, because $A$
is abelian and $P$ is a
permutation module, hence $k(P)^A$ is pure over $k$ by Fischer's theorem.
\end{proof}

\begin{prop}                 \label{prop.multiplicative.bru.formula}
Let $A$ be a bicyclic group, and let
$M$ be an $A$-lattice.
Then $\Bru k(M)^A =0$.
\end{prop}

\begin{proof}
The proof of this proposition is much more
technical. Indeed it is the core of the proof of Theorem \ref{thm.multiplicative.bru.formula}. We shall refer
to Saltman's original
paper \cite{Sa10} and to Barge's more natural proof \cite{Ba1}.
\end{proof}

\begin{rem}                 \label{rem.multiplicative.bru.formula}
That the kernel
$\ker[H^2(G,k^*\oplus M)\xrightarrow{\Res} \prod_{A\in\C{B}_G} H^2(A,k^*\oplus
M)] $ is a subgroup of $\Bru k(M)^G$ does not rely on the last two
propositions. This accounts for
Saltman's early counter-examples to the ``Noether problem'' over an
algebraically closed field \cites{Sa3,Sa6}.
\end{rem}

\begin{prop}[Saltman \cite{Sa6}]                   \label{prop.multiplicative.bru.formula.example}
Let $G$ be a finite
group of order $p^n$, $p$
prime, and assume that $G$ is neither cyclic nor bicyclic -- which implies
$n\geq 3$. Let $M$ be a faithful
$G$-lattice and assume $\exp H^2(G,M)=p^n$. Then $\Bru k(M)^G\neq 0$.
\end{prop}

\begin{proof}
By hypothesis there exists an
element $\alpha\in
H^2(G,M)$ of order $p^n$. Then
\[
0\neq\beta:=p^{n-1}\alpha\in H^2(G,M).
\]
If $A$ is a cyclic or a bicyclic subgroup subgroup of $G$, then $A\neq G$ hence
$\#A=p^m$ with
$m\leq n-1$. In particular, the group $H^2(A,M)$ is killed by $p^{n-1}$ and the
restriction of
$\beta=p^{n-1}\alpha$ to $H^2(A,M)$ is trivial. Hence
\[
0\neq\beta\in\ker\Bigl\lbrack H^2(G,M)\xrightarrow{\Res}\prod_{A\in\C{B}_G}H^2(A,M)\Bigr\rbrack\subset
\Bru k(M)^G
\]
where the last inclusion is part of Theorem \ref{thm.multiplicative.bru.formula}.
\end{proof}

\begin{exa}[multiplicative  example \cite{Sa6}]                 \label{exa.multiplicative.bru.formula}
Let $G$ be as above,
i.e.\ of order
$\#G=p^n$, $p$ prime and $G$ not bicyclic which implies $n \geq 3$. In the
standard
resolution of the trivial $G$-lattice $\D{Z}$, let $M$ denote the kernel of
the map $\pi\colon\D{Z}\croch{G\times G}\to\D{Z}\croch{G}$ given on generators of $\D{Z}\croch{G\times G}$ by $(g,h)\mapsto
(g-h)$. From the exact
sequence 
\[
0 \to  M\xrightarrow{}\D{Z}\croch{G\times G}\xrightarrow{\pi}\D{Z}\croch{G}\xrightarrow{}
\D{Z} \to  0
\]
we deduce $H^2(G,M) \cong\widehat H^0(G,\D{Z}) \cong \D{Z}/p^n$ where $\#G=p^n$.
Then 
\[
\exp H^2(G,M)=p^n. 
\]
This gives for every such group $G$ a faithful $G$-lattice $M$
such that
\[
\Bru k(M)^G \neq 0.
\]
For a given prime $p$ the simplest group is $G=(\D{Z}/p)^3$.
\end{exa}

\begin{exa}[linear example]                    \label{exa.linear.bru.formula}
Using Lemma \ref{lem.multiplicative.action.linear},
for any prime $p$
this example leads to a linear
example for a group $\Gamma$
of order $p^{3p^3}$.
\end{exa}

\begin{cor}[Barge \cite{Ba1}]                 \label{cor.barge.multiplicative.bru.formula}
For a finite group $G$,
the following conditions
are equivalent:
\begin{enumerate}[\rm(i)]
\item
For every faithful
$G$-lattice $M$, we have $\Bru
k(M)^G=0$.
\item
All Sylow subgroups of $G$ are bicyclic.
\end{enumerate}
\end{cor}

\begin{proof}
Let us prove that (ii)$\implies$(i). Let $p$ be a prime and $G_p\subset G$ a $p$-Sylow subgroup.
The field extension $k(M)^{G_p}/k(M)^G$ is of degree prime to
$p$. It therefore induces an imbedding on $p$-primary components of Brauer groups, and also of unramified Brauer groups:
\[
\Bru(k(M)^G)(p\subset\Bru(k(M)^{G_p})(p,
\]
and the later group vanishes according to Proposition \ref{prop.multiplicative.bru.formula} applied to the bicyclic
group $G_p$.

For the converse, assume there exists an $\ell$-Sylow subgroup $G_{\ell}\subset G$ which is not bicyclic.
Let $\ell^n$ be its order. Let $M$ be the faithful $G$-lattice defined by the exact sequence from Example \ref{exa.multiplicative.bru.formula}:
\[
0 \to  M\xrightarrow{}\D{Z}\croch{G\times G}\xrightarrow{\pi}\D{Z}\croch{G}\xrightarrow{}
\D{Z} \to  0.
\]
Then 
\[
H^2(G,M)\cong\oplus_p\D{Z}/p^{n_p}
\]
where $\prod_pp^{n_p}$ is the order of the group $G$. Let $\beta\in H^2(G,M)$ be the element defined 
in $\oplus_p\D{Z}/p^{n_p}$ by
\[
\beta_p=\begin{cases}
                \ell^{n-1}&\text{if }p=\ell,\\
                0&\text{if }p\neq\ell.
                \end{cases}
\]
The same argument as in Proposition \ref{prop.multiplicative.bru.formula.example} shows that $\beta\neq 0$
is unramified, hence $\Bru k(M)^G\neq 0$.
\end{proof}

\begin{thm}[Saltman \cite{Sa10}]                     \label{thm.twisted.multiplicative.bru.formula}
Let $G$ be a finite
group, let $M$ be a faithful
$G$-lattice and $\alpha\in Z^1(G,\widehat M)$. Let $k(M)_\alpha$ denote the
field of fractions of
the twisted $G$-$k$-algebra $k\croch{M}$. One then has:
\[
\Bru(k(M)_\alpha^G)\cong	\ker\Bigl\lbrack H^2(G,k\croch{M}_\alpha^*) \xrightarrow{\Res}
\prod_{A\in\C{B}_G'}
H^2(A,k^*\oplus M)\Bigr\rbrack
\]
where $\C{B}_G'$ denotes the set of bicyclic subgroups $A$ of $G$ such that the
map
\[
H^2(A,k^*)\xrightarrow{} H^2(A,k\croch{M}_\alpha^*)
\]
is injective.
\end{thm}

\begin{proof}
Once more we refer to Saltman's original
paper \cite{Sa10} and to Barge's more natural proof \cite{Ba2}.
\end{proof}

\begin{thm}[Barge \cite{Ba2}]       \label{thm.twisted.multiplicative.bru.triviality}
For a finite group $G$,
the following conditions
are equivalent:
\begin{enumerate}[\rm(i)]
\item
For every faithful $G$-twisted multiplicative action $(M,\alpha)$, we have
\[
\Bru k(M)_\alpha^G=0.
\]
\item
All Sylow subgroups of $G$ are cyclic.
\end{enumerate}
\end{thm}

\begin{exa}[twisted multiplicative example]      \label{exa.twisted.multiplicative.obstruction}
Let $p$ be a prime and let 
$G\cong\D{Z}/p\times\D{Z}/p$. There exists a   twisted multiplicative   action of $G$ on a lattice 
$M$ such that
$\Bru k(M)_\alpha^G \neq 0$.
 In particular  there exists a twisted multiplicative action of  $G=\D{Z}/2\times\D{Z}/2$
on a lattice $M$ such that the field of invariants $k(M)_\alpha^G$ is not pure.
\end{exa}

\begin{rem}                                                       \label{rem.twisted.multiplicative.action}
Given a $G$-lattice $M$ one may show that the field
$k(M)^G$ is 
stably equivalent to the function field of a torus defined on the field $k(\D{Z}[G])^G$.
\end{rem}

For some computations on multiplicative invariants, see \cites{Be1,Be2,Be3,Be4,HaKan1,HaKan2,HaKan4,Sa9}.
For a recent general report on multiplicative actions, we refer to the forthcoming book by M.~Lorenz \cite{Lor}.

\section{Homogeneous spaces}
                        \label{sec.homogeneous.spaces}

In this section, the ground field $k$ is the field $\D{C}$ of complex numbers.
We fix an isomorphism between the group of all roots of unity in $\D{C}$
and $\D{Q}/\D{Z}$.
The tools used are topological. However standard arguments
show that once the statements of Theorems \ref{thm.almost.free.bru.triviality} and
\ref{thm.homogeneous.space.bru.triviality} have been proved over
$\D{C}$, they hold over any algebraically closed field $\Omega$ of characteristic zero:
reduction to an algebraically closed field $k$ which may be embedded both in $\D{C}$ and $\Omega$,
then invariance of \'etale cohomology with finite coefficients under extension of algebraically closed ground fields.

\subsection{The case of an almost free linear 
                                                      representation}
                        \label{sec.almost.free.representation}

\begin{thm}[Bogomolov \cite{Bo3}*{Lemma 5.7}]                         \label{thm.almost.free.bru.triviality}
Let $G$ be a connected
algebraic group over the complex field $k=\D{C}$, and let $V$ be an
almost free finite dimensional linear
representation of $G$. Then $\Bru k(V)^G=0$.
\end{thm}

For   $G=\PGL_n$, this is a theorem of Saltman  \cite{Sa5}*{Theorem 2.9}, for which alternate proofs are given in \cite{CS2}*{Theorem 9.7} and \cite{Sa11}. This answered a question of Procesi \cite{Pro2}.

\begin{lem}                                                                                                   \label{lem.very.good.representation}
For any linear algebraic group $G$ over a field $k$ and 
any positive integer $s$, there exists a 
$k$-linear representation $V$ of $G$ and a non-empty $G$-stable
open set $U\subset V$ such that
\begin{enumerate}[\rm(i)]
\item
the complement $V\setminus U$ is of codimension $\geq s$,
\item
there is a
morphism $U\to U/G$ of $k$-varieties making $U$ 
into a $G$-torsor over $U/G$.
\end{enumerate}
\end{lem}

\begin{proof}
See Totaro \cite{To}*{Remark 1.4, p.~252}. Let $G\subset\GL_n$ be any faithful $k$-linear representation of $G$.
Let $N\geq 1$ be an integer. Let $W=\C{M}_{N+n}$ be the vector space of $(N+n)\times(N+n)$ matrices $M$ over $k$:
\[
M=\left(\hskip -.5\arraycolsep
     \begin{array}{c@{}c}
     \begin{array}{|ccc|}\hline
     &&\\[-2mm]
     &A&\\[-2mm]
     &&
     \end{array}&
     \begin{array}{c|}\hline
     \phantom{B}\\[-2mm]
     B\\[-2mm]
     \phantom{B}
     \end{array}\\
     \begin{array}{|ccc|}\hline
     &C&\\
     \hline
     \end{array}&
     \begin{array}{c|}\hline
     D\\
     \hline
     \end{array}
     \end{array}
     \hskip -.5\arraycolsep\right) \hskip -1.8\arraycolsep
\begin{array}{l@{}c}
\biggm\updownarrow&N\\[2.5mm]
\bigm\updownarrow&n
\end{array}
\]
Let $V=\C{M}_{n,N+n}$ be the vector space of $n\times(N+n)$ matrices over $k$ and let
$\pi\colon W\to V$ be the linear projection given by
\[
\left(\hskip -.5\arraycolsep
     \begin{array}{c@{}c}
     \begin{array}{|ccc|}\hline
     &&\\[-2mm]
     &A&\\[-2mm]
     &&
     \end{array}&
     \begin{array}{c|}\hline
     \phantom{B}\\[-2mm]
     B\\[-2mm]
     \phantom{B}
     \end{array}\\
     \begin{array}{|ccc|}\hline
     &C&\\
     \hline
     \end{array}&
     \begin{array}{c|}\hline
     D\\
     \hline
     \end{array}
     \end{array}
     \hskip -.5\arraycolsep\right)
\longmapsto
\left(\hskip -.5\arraycolsep
     \begin{array}{c@{}c}
     \begin{array}{|ccc|}\hline
     &C&\\
     \hline
     \end{array}&
     \begin{array}{c|}\hline
     D\\
     \hline
     \end{array}
     \end{array}
     \hskip -.5\arraycolsep\right).
\] 
Let $\Omega=\GL_{N+n}\subset W$ be the group of invertible matrices.
We denote by $\widetilde G\cong G$ the subgroup of matrices of the form
\[
\left(\hskip -.5\arraycolsep
     \begin{array}{c@{}c}
     \begin{array}{|ccc|}\hline
     &&\\[-2mm]
     &I&\\[-2mm]
     &&
     \end{array}&
     \begin{array}{c|}\hline
     \phantom{B}\\[-2mm]
     0\\[-2mm]
     \phantom{B}
     \end{array}\\
     \begin{array}{|ccc|}\hline
     &0&\\
     \hline
     \end{array}&
     \begin{array}{c|}\hline
     D\\
     \hline
     \end{array}
     \end{array}
     \hskip -.5\arraycolsep\right)
\]
where $D\in G$ and $I=I_N$ is the $N\times N$ identity matrix. We denote by
$H\vartriangleleft\GL_{N+n}$ the invariant subgroup consisting of matrices of the form
\[
\left(\hskip -.5\arraycolsep
     \begin{array}{c@{}c}
     \begin{array}{|ccc|}\hline
     &&\\[-2mm]
     &A&\\[-2mm]
     &&
     \end{array}&
     \begin{array}{c|}\hline
     \phantom{B}\\[-2mm]
     B\\[-2mm]
     \phantom{B}
     \end{array}\\
     \begin{array}{|ccc|}\hline
     &\,0\,&\\
     \hline
     \end{array}&
     \begin{array}{c|}\hline
     \,I\,\\
     \hline
     \end{array}
     \end{array}
     \hskip -.5\arraycolsep\right)
\]
where $A\in\GL_N$ and $I=I_n$ is the $n\times n$ identity matrix. 
The group $H$ is clearly an extension of $\GL_N$ by 
a unipotent invariant subgroup whose underlying $k$-variety is an affine space $\D{A}^{Nn}$.
The product $\Gamma=H\widetilde G\cong H\rtimes G$
is a subgroup of $\GL_{N+n}$.
Let $U\subset V$ denote the open dense subset consisting of matrices of rank $n$.
Let $Z$ be its complement. 
The projection $\pi\colon W\to V$ induces a surjective morphism of $G$-$k$-varieties
\[
\pi_{\Omega}\colon\Omega\to U.
\]
This map induces a $G$-isomorphism of $k$-varieties $\GL_{N+n}/H\cong U$. 
Since $\Gamma=H\rtimes G$ is a closed subgroup of $\GL_{N+n}$, the canonical morphism
\[
\GL_{N+n}\to\GL_{N+n}/\Gamma
\] 
factorizes through $\pi_{\Omega}$ and induces a morphism $\varpi$ of $k$-varieties
which makes $\GL_{N+n}/H$ into a $G$-torsor over $\GL_{N+n}/\Gamma$.
Thus $U/G$ exists and $U$ is a $G$-torsor over $U/G$:
\[
\begin{CD}
         \GL_{N+n}/H    @>\varpi>>  \GL_{N+n}/\Gamma\\
         @V\cong V\bar\pi_{\Omega}V                                 @V\cong VV\\
          U@>>> U/G.
\end{CD}
\]
Since $\codim_VZ$ goes to infinity as $N$
goes to infinity, we obtain $\codim_VZ\geq s$ for $N$ large.
\end{proof}

\begin{proof}[Proof of Theorem \ref{thm.almost.free.bru.triviality}]
In view of the no-name lemma and Proposition \ref{prop.bru.purity.invariance}, we may replace the linear representation
$V$ by a ``better one''. Namely, we take $V$ as in Lemma \ref{lem.very.good.representation}
for $s=3$.
Thus
we may now assume that we
are given a linear representation $V$ of $G$ and a $G$-stable open set $U$
with complement $Z$ such
that $\codim_VZ\geq 3$ and such that there is a morphism $U\to U/G$ which
makes $U$ into a $G$-torsor over $U/G$.

The connected group $G$ is an extension
\[
1\to G'\to G\to T\to 1
\]
of a torus $T$ by a connected group $G'$ without characters: $\Char(G')=0$. Now we may first
quotient $U$ by the action of $G'$,
and then by the action of $T$, getting a morphism $U/G'\to U/G$
which one checks 
makes
$U/G'$ into a $T$-torsor over $U/G$. Since torsors under tori are
locally trivial, $U/G'$ is birational to $T\times U/G$. The stability of the unramified Brauer
group under pure
extensions (Proposition \ref{prop.bru.purity.invariance}) now implies $\Bru (k(V)^G)\cong \Bru(k(V)^{G'})$. It is
thus enough to prove the
theorem when $G$ is connected and satisfies $\Char(G)=0$, which we now assume.

\begin{claim*} 
For such $G$ we now claim that
\begin{equation}                                                                  \label{eq.claim.formula}
\Bru k(V)^G \cong\ker\Bigl\lbrack H^2(BG,\D{Q}/\D{Z}) \to\prod_{A\in\C{B}_G}
H^2(BA,\D{Q}/\D{Z})\Bigr\rbrack.
\end{equation}
Here $BG$, resp.\ $BA$, denotes
the classifying space of the topological group $G=G(\D{C})$, resp.\ of the
finite group $A$. Recall that $\C{B}_G$ denotes the set of finite bicyclic subgroups of $G(\D{C})$.
\end{claim*}

Theorem \ref{thm.bru.bicyclic.reductive} 
implies
\[
\Bru (k(V)^G) = \accol{\alpha\in\Br k(V)^G\mid\alpha_A\in\Bru (k(V)^A) \text{ for all } A\in\C{B}_G}.
\]
Since $A$ is abelian and the action on $V$ is linear, $k(V)^A$ is pure by
Fischer's Theorem (Proposition \ref{prop.torus.purity}),
hence $\Bru k(V)^A=0$ (Proposition \ref{prop.bru.purity.invariance}). 
The field $k(V)^G$, resp.~$k(V)^A$, is the function field of the smooth variety $U/G$,
resp.~ $U/A$.
Using Theorem \ref{thm.groth.brauer.unramified}, we get:
\begin{equation}                                                                  \label{eq.open.brauer.nr.formula}
\Bru k(V)^G \cong	
\ker\Bigl\lbrack\Br(U/G) \to\prod_{A\in\C{B}_G} \Br(U/A)\Bigr\rbrack.
\end{equation}
\renewcommand{\qed}{}
\end{proof}

\begin{lem}                                                   \label{lem.discrete.brauer.group}
Let $X$ be a smooth algebraic variety
over $\D{C}$. If
$\Pic X$ is torsion, then we have a canonical isomorphism
\[
H_{\textup{\'et}}^2(X,\D{Q}/\D{Z}) \cong \Br X.
\]
\end{lem}

\begin{proof} 
This is a consequence of ``Kummer's exact sequence'':
\[
0 \to\Pic X\otimes \D{Q}/\D{Z} \to  H_{\text{\'et}}^2(X,\D{Q}/\D{Z}) \to\Br X \to 0.\qedhere
\]
\end{proof}

\begin{lem}                                                    \label{lem.pic.exact.sequence}
Let $X$ be a smooth connected variety over $k=\D{C}$.
Let $G$ be a linear algebraic group with character group $\Char(G)=\Hom(G,\D{G}_m)$.
Let $X\to X/G$ be a $G$-torsor.
\begin{enumerate}[\rm(i)]
\item
If $G$ is connected, there is a natural exact sequence
\[
1 \to  k[X/G]^* \to  k[X]^* \to  \Char(G) \to\Pic(X/G) \to\Pic X.
\]
\item
If $G$ is a finite constant group $A$, there is a natural exact sequence
\[
1\to\Char(A) \to\Pic(X/A) \to\Pic X.
\]
\end{enumerate}
\end{lem}

\begin{proof} 
For (i), see \cite{San}*{Proposition 6.10}. For (ii), the same arguments yield an exact sequence
\[
1\to\Hom(A,k\croch{X}^*)\to\Pic(X/A) \to\Pic X.
\]
Then use the fact that $k\croch{X}^*/k^*$ is torsionfree.
\end{proof}

\begin{lem}                                                       \label{lem.brauer.semisimple.quotient}
Let $X$ be a smooth connected variety over $\D{C}$.
Let $G$ be a linear algebraic group. 
Let $X\to X/G$ be a $G$-torsor.
Assume that $\Pic X$ is a torsion group and that $G$ is either finite, or connected and characterfree.
Then  we have a canonical isomorphism
\[
\Br(X/G) \cong H_{\textup{\'et}}^2(X/G,\D{Q}/\D{Z}).
\]
\end{lem}

\begin{proof} 
Under our assumptions, Lemma \ref{lem.pic.exact.sequence}
implies that $\Pic(X/G)$ is torsion. Then Lemma \ref{lem.discrete.brauer.group} gives the
isomorphism $\Br(U/G) \cong
H^2(U/G,\D{Q}/\D{Z})$.
\end{proof}

\begin{proof}[Proof of Theorem \ref{thm.almost.free.bru.triviality} (continued)]
Lemma \ref{lem.brauer.semisimple.quotient} applied to $X=U$ together with \eqref{eq.open.brauer.nr.formula}
implies
\[
\Bru k(V)^G \cong	\ker\Bigl\lbrack H^2(U/G,\D{Q}/\D{Z}) \to\prod_{A\in\C{B}_G}
H^2(U/A,\D{Q}/\D{Z})\Bigr\rbrack.
\]
Since $U\to U/G$ is a $G$-torsor, there is a map of topological
spaces $U/G\to BG$ 
such that $U\to U/G$ is the pull-back of the universal
covering space $EG\to BG$.
Similarly, there is a map $U/A\to BA$. Now the assumption on the
codimension of the complement of
$U$ in $V$, which had not been used yet, implies that the maps $U/G\to BG$
and $U/A\to BA$ induce
isomorphisms on the cohomology groups $H_{\text{\'et}}^i(\;\cdot\;,\D{Q}/\D{Z})$ for $i\leq 2$ (the
projections $U\to U/G$ and
$U\to U/A$ are algebraic approximations of the universal covering spaces).
The ``claim'' \eqref{eq.claim.formula} now follows.
\renewcommand{\qed}{}
\end{proof}

\begin{lem}                                    \label{lem.twisted.multiplicative.bru.triviality}
If $\widetilde G$ is a connected, simply
connected group, then
\[
H^1(B\widetilde G,\D{Q}/\D{Z})=H^2(B\widetilde G,\D{Q}/\D{Z})=0.
\]
\end{lem}

\begin{proof} 
The
fibration $E\widetilde G\to B\widetilde G$ yields isomorphisms 
$ \pi_1B\widetilde G=\pi_0\widetilde G = 0$ (since $\widetilde G$ is connected) and 
$\pi_2B\widetilde G=\pi_1\widetilde G=0$ (since $\widetilde G$ is simply connected).
Thus $B\widetilde G$ is $2$-connected. The Hurewicz theorem now yields
$H_1(B\widetilde G)=H_2(B\widetilde G)=0$. The universal coefficient theorem now
yields the statement of the lemma.
\end{proof}

\begin{lem}                                                                         \label{lem.htwo.classifiant}
If $G$ is a connected  group without characters, the
universal covering
$\widetilde G \to G$ defines a natural isomorphism of finite abelian groups
\[
d_{2,G}\colon\Hom(\pi_1G,\D{Q}/\D{Z}) \xrightarrow{\cong} H^2(BG,\D{Q}/\D{Z}).
\]
\end{lem}

\begin{proof} 
Since $G$ is a connected  group without characters, $G=R_{\text{u}}(G)\rtimes G_{\text{ss}}$ with $R_{\text{u}}(G)$ unipotent and $G_{\text{ss}}$ semisimple. Thus its fundamental group $\pi_1G$ 
is a finite abelian group.
We have the Leray-Serre spectral sequence \cite{Mcl}*{\S5, Theorem 5.2} for the fibration $\pi_1G\to B\widetilde G\to BG$:
\[
E_2^{p,q}=H^p(BG,\C{H}^q(\pi_1G,\D{Q}/\D{Z})) \implies H^{p+q}(B\widetilde G,\D{Q}/\D{Z}).
\]
Since $BG$ is $1$-connected, each local system $\C{H}^q(\pi_1G,\D{Q}/\D{Z})$ is constant. Then
the associated exact sequence of terms of lower degree
\[
H^1(B\widetilde G,\D{Q}/\D{Z})\to H^1(\pi_1G,\D{Q}/\D{Z}) \xrightarrow{d_{2,G}}
H^2(BG,\D{Q}/\D{Z})\to
H^2(B\widetilde G,\D{Q}/\D{Z})
\]
together with Lemma \ref{lem.twisted.multiplicative.bru.triviality} yields 
a natural isomorphism of finite abelian groups
$d_{2,G}\colon\Hom(\pi_1G,\D{Q}/\D{Z}) \xrightarrow{\cong} H^2(BG,\D{Q}/\D{Z})$.
\end{proof}

\begin{lem}                                                                         \label{lem.htwo.classifiant.fini}
Let $A$ be a finite group and $C$ be a trivial $A$-module.
Let
\begin{equation}                                                     \label{eq.central.extension}
1\to\pi\to\widetilde A\to A\to 1\tag{$\C{E}_A$}
\end{equation}
be a central extension. The associated Lyndon-Hochschild-Serre spectral sequence 
defines a natural morphism
\[
d_{2,\C{E}_A}\colon\Hom(\pi,C)\to H^2(A,C)
\]
sending $\phi\colon\pi\to C$ to the class of the extension $\phi_*(\C{E}_A)$ obtained from  \eqref{eq.central.extension}
by push-out along $\phi$.
\end{lem}

\begin{proof}
The Leray-Serre spectral sequence 
\cite{Mcl}*{\S5, Theorem 5.2} for the fibration $\pi\to B\widetilde A\to BA$ associated to \eqref{eq.central.extension}:
\[
E_2^{p,q}=H^p(BA,\C{H}^q(\pi,C))\implies H^{p+q}(B\widetilde A,C)
\]
coincides 
with the Lyndon-Hochschild-Serre spectral sequence (see \cite{Mcl}*{\S$8^{\textit{bis}}.2$, p.~ 342})
\[
E_2^{p,q}=H^p(A,H^q(\pi,C))\implies H^{p+q}(\widetilde A,C)
\]
associated to \eqref{eq.central.extension}. This extension being central, the action of $A$ on $H^q(\pi,C)$ is trivial and the local system $\C{H}^q(\pi,C)$ is constant on $BA$. Hence both morphisms
$H^0(BA,\C{H}^1(\pi,C))\to H^2(BA,\C{H}^0(\pi,C))$ and $H^0(A,H^1(\pi,C))\to H^2(A,H^0(\pi,C))$ coincide,  
yielding the  natural morphism
\[
d_{2,\C{E}_A}\colon\Hom(\pi,C)\to H^2(A,C).
\]
It is well known (see \cite{Bro}*{\S IV.3, Theorem 3.12}) that 
$H^2(A,C)$ classifies abstract central extensions of $A$ by $C$. Let $C=\pi$. One can check that 
$d_{2,\C{E}_A}(\id_{\pi})$ is the class of the central extension $\C{E}_A$. 
Now for any $C$, the functoriality of the LHS spectral sequence implies that $d_{2,\C{E}_A}(\phi)$ is the class of $\phi_*(\C{E}_A)$.
\end{proof}

\begin{proof}[Proof of Theorem \ref{thm.almost.free.bru.triviality} (continued)]
Let $R_{\text{u}}(G)$ be the unipotent radical of $G$. 
We have the exact sequence
\begin{equation}                                                                                             \label{eq.reduct.semisimple}
1 \to R_{\text{u}}(G)\to G\xrightarrow{\varrho}G_1\to 1,
\end{equation}
where $G_1$ is semisimple.
Pulling back the universal cover $\widetilde G_1\to G_1$ over $G$, we get the commutative diagram 
of exact sequences
\[
\begin{CD}
         1         @>>>  R_{\text{u}}(G) @>>>    \widetilde G     @>>>   \widetilde G_1            @>>>   1\\
         @.                    @|                          @VVV                        @VVV                @.  \\
         1         @>>> R_{\text{u}}(G)  @>>>      G    @>>>   G_1              @>>>    1,
\end{CD}
\]
where $\widetilde G\to G$ is the universal cover of $G$.
From Lemma \ref{lem.htwo.classifiant} we get the commutative diagram
\[
\begin{CD}
         \Hom(\pi_1G,\D{Q}/\D{Z})@>\cong>>   H^2(BG,\D{Q}/\D{Z})\\
         @AAA                                 @AAA \\
          \Hom(\pi_1G_1,\D{Q}/\D{Z})@>\cong>>   H^2(BG_1,\D{Q}/\D{Z}).
\end{CD}
\]
Since $R_{\text{u}}(G)$ is contractible, the map $\pi_1G\to\pi_1G_1$ induced by \eqref{eq.reduct.semisimple}
is an isomorphism. Hence the map $H^2(BG_1,\D{Q}/\D{Z})\to H^2(BG,\D{Q}/\D{Z})$ is an isomorphism.
For each finite subgroup $A\subset G$, the exact sequence \eqref{eq.reduct.semisimple} induces an isomorphism $A\cong\varrho(A)$. Moreover it induces a bijection $A\mapsto\varrho(A)$
between the finite subgroups of $G$ and those of $G_1$. We conclude that
\[
\ker\Bigl\lbrack H^2(BG,\D{Q}/\D{Z}) \to\prod_{A\in\C{B}_G}H^2(BA,\D{Q}/\D{Z})\Bigr\rbrack
\]
and
\[
\ker\Bigl\lbrack H^2(BG_1,\D{Q}/\D{Z}) \to\prod_{A\in\C{B}_{G_1}}H^2(BA,\D{Q}/\D{Z})\Bigr\rbrack
\]
are isomorphic. To prove that the first  kernel is trivial we may therefore assume that $G$ is a semisimple group.

Consider now the natural isogeny $\widetilde G\to G$ where $\widetilde G$
denotes
the simply connected cover of $G$. Each $A\in\C{B}_G$ gives rise to a
commutative diagram of exact sequences of groups
\[
\begin{CD}
         1         @>>>   \pi_1G @>>>    \widetilde G     @>>>   G              @>>>   1\\
         @.                    @|                          @AAA                        @AAA                @.  \\
         1         @>>>  \pi_1G  @>>>      \widetilde A    @>>>   A              @>>>    1.
\end{CD}
\]
and to a commutative square of fibrations
\begin{equation}                                                       \label{eq.fibrations.square}
\begin{CD}
         B\widetilde G     @>>>   BG\\
         @AAA                              @AAA           \\
         B\widetilde A    @>>>   BA.
\end{CD}
\end{equation}
Then \eqref{eq.fibrations.square}  
gives a morphism of Leray-Serre spectral sequences:
\[
\begin{array}{ccc}
E_2^{p,q}=H^p(BG,\C{H}^q(\pi_1G,\D{Q}/\D{Z})) & \implies & H^{p+q}(B\widetilde G,\D{Q}/\D{Z}) \\
\downarrow &&  \downarrow\\
E_2^{p,q}=H^p(BA,\C{H}^q(\pi_1G,\D{Q}/\D{Z})) & \implies & H^{p+q}(B\widetilde A,\D{Q}/\D{Z}),
\end{array}
\]
hence $d_2\colon E_2^{0,1}\to E_2^{2,0}$ 
defines a commutative square (see Lemmas  \ref{lem.htwo.classifiant} and \ref{lem.htwo.classifiant.fini})
\begin{equation}                                                                   \label{eq.square.htwo.classifiant}
\begin{CD}
\Hom(\pi_1G,\D{Q}/\D{Z})  @>d_{2,G}>>   H^2(BG,\D{Q}/\D{Z})\\
    @|                                                                   @VVV                    \\
\Hom(\pi_1G,\D{Q}/\D{Z})  @>d_{2,\C{E}_A}>>  H^2(BA,\D{Q}/\D{Z}).
\end{CD}
\end{equation}
Any element $\alpha\in H^2(BG,\D{Q}/\D{Z})$ may be
interpreted as a homomorphism $\phi$ of $\pi_1G =\pi_1^{\text{alg}}G$ to $\D{Q}/\D{Z}$,
hence to some
$\D{Z}/n=\im\phi$. Hence by push-out from the natural extension  
\[
1\to
\pi_1G\to {\widetilde G}\to G\to 1
\]
it defines an extension 
\[
1\to\D{Z}/n \to  G_1 \to  G \to  1,
\]
where $G_1$ is a connected semisimple algebraic group over $\D{C}$. This extension
itself defines an
extension of $G$ by $\D{Q}/\D{Z}$.
By Lemma \ref{lem.htwo.classifiant.fini}, the natural restriction map $r\colon
H^2(BG,\D{Q}/\D{Z})\to
H^2(BA,\D{Q}/\D{Z})$ may be interpreted in terms of restrictions of extensions.
More precisely, if
$\alpha\in H^2(BG,\D{Q}/\D{Z})$ corresponds to the isogeny
\[
1\to\D{Z}/n \to  G_1 \to  G \to  1
\]
(notations as above), and if $A\subset
G$ is a finite subgroup, then the restriction of the above isogeny to this
subgroup gives rise to a
central extension
\[
1 \to\D{Z}/n\to  A_1 \to  A \to  1
\]
whose class in $H^2(A,\D{Z}/n)$ restricts to the
class of $r(\alpha)\in H^2(A,\D{Q}/\D{Z})\cong H^2(BA,\D{Q}/\D{Z})$.

To prove that
\[
\ker\Bigl\lbrack H^2(BG,\D{Q}/\D{Z}) \to\prod_{A\in\C{B}_G}H^2(BA,\D{Q}/\D{Z})\Bigr\rbrack=0,
\]
which will complete the proof of the theorem, all we now need to prove is that if we are given
a nontrivial central extension
\begin{equation}                                      \label{eq.nontrivial.central.extension}
1\to\D{Z}/n\to  G_1 \to  G \to  1
\end{equation}
with $G_1$ connected, then there exists a finite bicyclic subgroup $A$ in the
semisimple group $G$ such that
the restriction
\begin{equation}                                      \label{eq.restriction.central.extension}
1\to\D{Z}/n \to  A_1 \to  A\to  1
\end{equation}
of \eqref{eq.nontrivial.central.extension} to $A$ is a nontrivial extension. 
Since \eqref{eq.restriction.central.extension} is a
central extension, it is nontrivial if $A_1$ is not commutative. 
\renewcommand{\qed}{}
\end{proof}

Let us postpone the proof of the:

\begin{prop}                                                     \label{prop.central.commutator}
Any
element in the centre of a connected semisimple group is a commutator of
elements of finite order.
\end{prop}

\begin{proof}[Proof of Theorem \ref{thm.almost.free.bru.triviality} (end)]
Let $c$ be a generator of the subgroup $\D{Z}/n$ of the centre of the connected semisimple
group $G_1$.
Using Proposition \ref{prop.central.commutator}
we can write $c$ as a commutator
$a_1b_1a_1^{-1}b_1^{-1}$, with $a_1$ and $b_1$ of finite order in $G_1$, let
$a$ and $b$ be the images
of $a_1$ and $b_1$ in $G$, and let $A\subset  G$ be the finite, abelian,
bicyclic group which they
generate. Since $a_1$ and $b_1$ do not commute, the group $A_1$ is not
commutative and the proof of
Theorem \ref{thm.almost.free.bru.triviality} is complete.
\end{proof}

\begin{proof}[Proof of Proposition \ref{prop.central.commutator}] 
It is
enough to give the proof when the connected semisimple group $G$ is simply
connected, and then it is
enough to prove it when $G$ is simple and simply connected.
\renewcommand{\qed}{}
\end{proof}

\begin{proof}[First proof (suggested by O. Gabber)] 
To
each choice of a basis of roots in the character group $\Char(T)$ of a maximal
torus $T\subset  G$ one
may associate an element $c$ of the Weyl group $W$ of $G$, known as the
Coxeter element \cite{Bbk}*{V.6.2 and VI.1.11}.

It is known
(loc. cit.) that $1$ is not an eigenvalue of $c$ for its action on
$\F{t}=\Char(T)\otimes\D{C}$. Let $d$ in the
normalizer $N(T)$ of $T$ be a representant of $c$. Such a representant may
be chosen of finite order.
Conjugacy by $d$ on $T$ induces an automorphism whose tangent linear map is
$c\in\Aut\F{t}$. Now
the tangent map to the homomorphism $\lambda\colon T \to  T$ given by $\lambda(t)=
dtd^{-1}t^{-1}$ is
$c-\id\in\End\F{t}$, hence is invertible. Thus $\lambda$ is an isogeny,
and this isogeny
induces a surjection on torsion points $\lambda\colon T_{\text{tors}}\to  T_{\text{tors}}$. We conclude that any
element $x\in T_{\text{tors}}$ may be written as $x=dyd^{-1}y^{-1}$ in $G$,
with $y\in T_{\text{tors}}$.
Since the centre of $G$ is contained in any maximal torus of $G$, the
conclusion follows.
\end{proof}

\begin{proof}[Second proof (suggested by J-P. Serre)] 
One first
proves the result for $G=\SL_n$ (a slick proof being as follows: the
algebra generated by $a$ and
$b$ with the relations $a^n=b^n=1$ and $ab=\zeta ba$ with $\zeta^n = 1$ is
none other than
$\C{M}_n(\D{C}))$. Now inspection of the root systems reveals that any simply
connected semisimple
group $G$ contains a subgroup $H\cong\SL_{n_1}\times\dots\times
\SL_{n_r}$ with $\rank H=\rank G$,
hence with $\centre(G)\subset\centre(H)$.
\end{proof}
 
\subsection{The case of a homogeneous space}
                                         \label{ssec.homogeneous space}
                        
For Theorem \ref{thm.homogeneous.space.bru.triviality}, we need some preparation.
Let us first recall a theorem of Steinberg
(cf.\ \cite{SpSt}*{II, 3.9, p.~197}).

\begin{thm}[Steinberg]                                \label{thm.steinberg}
If $G$ is a semisimple simply connected
group, the centralizer of a
semisimple element of $G$ is a connected reductive group.
\end{thm}

\begin{rem*}
This is no longer true if $G$ is not simply connected.
\end{rem*}

\begin{cor}                                         \label{cor.steinberg}
Given two commuting semisimple
elements in a semisimple, simply
connected group $G$, there is a maximal torus of $G$ which contains them
both.
\end{cor}

\begin{proof}
Let $x,y$ be two such elements. Let $H=Z_G(x)$ be the centralizer of $x$. Then $H$ is a connected reductive group (Theorem \ref{thm.steinberg}), and $x$ is in its center. Hence $x$ belongs to any maximal torus of $H$.
The element $y$ belongs to $H$. Since it is semisimple, it belongs to a maximal torus $T$ of $H$.
There exists a maximal torus of $G$ which contains $T$. Such a torus contains both $x$ and $y$.
\end{proof}

\begin{prop}                                                           \label{prop.bicyclic.quotient.bru.triviality}
Let $A$ be a bicyclic finite subgroup of a
simply
connected group $G$. Then $G/A$ is a rational variety, and the unramified
Brauer group of $G/A$ is
trivial.
\end{prop}

\begin{proof} 
We have the natural exact sequence
\[
1\to U\to G\xrightarrow{\pi}G_{\text{ss}}\to 1,
\]
where $U$ is the unipotent radical of $G$ and $G_{\text{ss}}$ is a connected, semisimple, simply connected group.
Since $U\cap A=\accol{1}$, the map $A\to\pi(A)$ is an isomorphism and the natural projection $G/A\to G_{\text{ss}}/A$
defines a $U$-torsor. Hence the varieties $G/A$ and $G_{\text{ss}}/A\times U$ are isomorphic. To prove the proposition
we may therefore assume that $G$ is semisimple and simply connected.

Let $A\subset  G$ be a
finite bicyclic group. By the previous corollary, there is a maximal torus
$T$ which contains $A$.
Let $T'=T/A$. The map $G\to G/T$ makes $G$ into a $T$-torsor over $G/T$,
hence $G$ is birational to
$T\times G/T$. On the other hand, the map $G/A\to G/T$ makes $G/A$ into a
$T'$-torsor over $G/T$,
hence $G/A$ is birational to $T'\times G/T$. Since $T$ and $T'$ are birational to
each other, we conclude
that $G/A$ is birational to $G$, hence is a rational variety.
\end{proof}

The proof of the following Theorem uses Theorem  \ref{thm.almost.free.bru.triviality}.
In turn, Proposition \ref{prop.special.stable.equivalence} shows that 
Theorem \ref{thm.homogeneous.space.bru.triviality} generalizes Theorem  \ref{thm.almost.free.bru.triviality}.

\begin{thm}[Bogomolov \cite{Bo4}*{Theorem 2.4}]                   \label{thm.homogeneous.space.bru.triviality}
Let $G$ be a connected, simply connected
group over
$\D{C}$, and let $H\subset  G$ be a connected closed subgroup. Then the
unramified Brauer group of
$G/H$ vanishes.
\end{thm}

\begin{proof}
There is a
closed normal subgroup $H_1\subset  H$ which is connected and
characterfree, such that
$H/H_1$ is a torus $T$. Thus the map $G/H_1\to G/H$ makes the
first variety into a
$T$-torsor over the second one. Since torsors under tori are locally
trivial, $G/H_1$ is
birational to $T\times G/H$. By Proposition \ref{prop.bru.purity.invariance} it is thus
enough to prove the  theorem when $H$ is connected and characterfree, which we now assume.

Just as in the proof of Theorem \ref{thm.almost.free.bru.triviality}, we have an
exact sequence:
\[
0\to\Pic(G/H)\otimes\D{Q}/\D{Z}\to  H_{\text{\'et}}^2(G/H,\D{Q}/\D{Z})\to\Br (G/H)\to 0
\]
and similar exact
sequences for any finite subgroup $A\subset H$:
\[
0\to\Pic(G/A)\otimes\D{Q}/\D{Z}\to  H_{\text{\'et}}^2(G/A,\D{Q}/\D{Z})\to\Br (G/A)\to 0.
\]
We also have exact sequences (Lemma \ref{lem.pic.exact.sequence})
\[
\Char(H)\to\Pic(G/H)\to\Pic G
\]
and
\[
\Char(A)\to\Pic(G/A)\to\Pic G.
\]
Now $\Pic G=0$ because $G$ is connected and simply connected, $\Char(H)=0$, and
$\Char(A)$ is finite since $A$ is finite. 
We thus get a commutative diagram:
\begin{equation}                                               \label{eq.square.discrete.brauer.group}
\begin{CD}
H_{\text{\'et}}^2(G/A,\D{Q}/\D{Z}) @> \cong>> \Br (G/A) \\
@AAA                                                 @AAA \\
H_{\text{\'et}}^2(G/H,\D{Q}/\D{Z}) @> \cong>> \Br (G/H).
\end{CD}
\end{equation}
Proposition \ref{prop.bicyclic.quotient.bru.triviality} together with the general formula for the unramified Brauer group (Theorem \ref{thm.bru.bicyclic.reductive}) 
and the above diagram imply that
\begin{equation}                                                      \label{eq.formula.brauernr.homogeneous.space}
\Bru G/H \cong	
\ker\Bigl\lbrack H^2(G/H,\D{Q}/\D{Z}) \to\prod_{A\in\C{B}_H}
H^2(G/A,\D{Q}/\D{Z})\Bigr\rbrack.
\end{equation}
From the universal property of the topological $H$-fibration
$EH\to BH$ there is a cartesian diagram of topological morphisms
\[
\begin{CD}
G @>>> EH \\
@VVV          @VVV\\
G/H @>>> BH.
\end{CD}
\]
This induces a map of spectral sequences
\[
\begin{array}{ccc}
E_2^{p,q}=H^p(G/H,\C{H}^q(H,\D{Q}/\D{Z})) & \implies & H^{p+q}(G,\D{Q}/\D{Z}) \\
\uparrow &&  \uparrow\\
E_2^{p,q}=H^p(BH,\C{H}^q(H,\D{Q}/\D{Z})) & \implies & H^{p+q}(EH,\D{Q}/\D{Z}),
\end{array}
\]
hence a commutative diagram:
\[
\begin{CD}
H^0(G/H,\C{H}^1(H,\D{Q}/\D{Z})) @> \partial>>  H^2(G/H,\C{H}^0(H,\D{Q}/\D{Z})) \\
@AAA                                                                  @AAA                   \\
H^0(BH,\C{H}^1(H,\D{Q}/\D{Z})) @> \partial>>  H^2(BH,\C{H}^0(H,\D{Q}/\D{Z})).
\end{CD}
\]
Since $G$ is simply connected, $\pi_1G=0$. It is also known (E.~Cartan's theorem) that $\pi_2G=0$. 
Thus $H^1(G,\D{Q}/\D{Z})=0$ and the universal
coefficient
theorem also gives $H^2(G,\D{Q}/\D{Z})=0$. We clearly have the same vanishing properties for $EH$, which 
is contractible.
This implies that the maps $\partial$ are isomorphisms.

Both $G/H$ and $BH$ are connected and simply connected, as follows
from the long sequence of
homotopy groups deduced from the two fibrations ($H$ is connected), hence
each local system $\C{H}^q(H,\D{Q}/\D{Z})$ is
constant and equal to $H^q(H,\D{Q}/\D{Z})$. Since $BH$ and $G/H$ are connected, the 
L.H.S.\ vertical map is identity on $H^1(H,\D{Q}/\D{Z})$, hence the R.H.S.\
vertical map is an isomorphism. But this map is the natural map $H^2(BH,\D{Q}/\D{Z})\to H^2(G/H,\D{Q}/\D{Z})$.

\begin{rem*} 
A shorter proof runs as follows. We may regard the map $G/H\to BH$ as a
topological fibration with fibre $G$. Indeed the map $G/H\times EG\to(EG)/H$ makes $G/H\times EG$ into a principal $G$-bundle over $(EG)/H$. 
The total space $G/H\times EG$ has the same homotopy type as $G/H$ and the base $(EG)/H$ has the same homotopy type as $BH$. Associated to this fibration there is a spectral sequence
\[
E_2^{p,q}=H^p(BH,H^q(G,\D{Q}/\D{Z})) \implies H^{p+q}(G/H,\D{Q}/\D{Z}).
\]
and this spectral sequence yields an isomorphism
\[
H^2(BH,\D{Q}/\D{Z}) \xrightarrow{\cong} H^2(G/H,\D{Q}/\D{Z})                                                                                                                                                                                                                                                
\]
which is induced by the
map $G/H\to BH$.
\end{rem*}

From the universal property of the topological $A$-covering
$EA\to BA$ there is a cartesian diagram of topological morphisms
\[
\begin{CD}
G @>>> EA \\
@VVV          @VVV\\
G/A @>>> BA.
\end{CD}
\]
This induces a map of spectral sequences
\[
\begin{array}{ccc}
E_2^{p,q}=H^p(A,H^q(G,\D{Q}/\D{Z})) & \implies & H^{p+q}(G/A,\D{Q}/\D{Z}) \\
\uparrow &&  \uparrow\\
E_2^{p,q}=H^p(A,H^q(EA,\D{Q}/\D{Z})) & \implies & H^{p+q}(BA,\D{Q}/\D{Z}).
\end{array}
\]
Since $H^1(G,\D{Q}/\D{Z})$ and $H^2(G,\D{Q}/\D{Z})$ are zero, and also all $H^q(EA,\D{Q}/\D{Z})$ for $q>0$,
we get the commutative diagram:
\[
\begin{CD}
H^2(A,H^0(G,\D{Q}/\D{Z})) @>\cong>>  H^2(G/A,\D{Q}/\D{Z})) \\
@AAA                                                                  @AAA                   \\
H^2(A,H^0(EA,\D{Q}/\D{Z})) @>\cong>>  H^2(BA,\D{Q}/\D{Z})).
\end{CD}
\]
Since $G$ and $EA$ are connected, the L.H.S.\ vertical map is identity on the group
$H^2(A,\D{Q}/\D{Z})$, hence the R.H.S.\ vertical map is an isomorphism:
\[
H^2(BA,\D{Q}/\D{Z}) \xrightarrow{\cong} H^2(G/A,\D{Q}/\D{Z}).
\]
Using the diagram
\[
\begin{CD}
G/A @>>> BA \\
@VVV          @VVV\\
G/H @>>> BH
\end{CD}
\]
we get the commutative diagram
\[
\begin{CD}
H^2(G/H,\D{Q}/\D{Z}) @>>>  H^2(G/A,\D{Q}/\D{Z})) \\
@AAA                                                                  @AAA                   \\
H^2(BH,\D{Q}/\D{Z}) @>>>  H^2(BA,\D{Q}/\D{Z})).
\end{CD}
\]
We have proved that both vertical maps are isomorphisms.
Formula \eqref{eq.formula.brauernr.homogeneous.space} now yields:
\[
\Bru G/H \cong	\ker\Bigl\lbrack H^2(BH,\D{Q}/\D{Z}) \to\prod_{A\in\C{B}_H}
H^2(BA,\D{Q}/\D{Z})\Bigr\rbrack,
\]
and we saw in the proof of Theorem \ref{thm.almost.free.bru.triviality}
that this kernel is zero.
\end{proof}

\begin{rem}                                            \label{rem.homogeneous.space.bru.triviality}
It would be interesting to know whether Theorem \ref{thm.homogeneous.space.bru.triviality}
extends to quotients $G/H$ with $G$ and $H$ connected. 

Note that if $G$ is connected, semisimple but not simply connected, and if $\widetilde G$ is its simply connected cover, the inverse image of
$H$ in $\widetilde G$ need not be connected, so that one may not reduce to the
situation of the
theorem. For example, if we consider in $\widetilde G=\SL(3)$ the subgroup $H_1$ which is the product of $\accol{1}\times\SL(2)$ 
by the diagonal $\mu_3$, the image $H$ of that group in the projective special linear group $\PSL(3)$ is clearly
connected whereas $H_1$ is not connected.
\end{rem}


\end{document}